\documentclass[11pt]{article}
\usepackage[utf8]{inputenc}

\usepackage[T1]{fontenc}

\usepackage{amsmath,amssymb,amsthm}
\usepackage{blindtext}
\usepackage{graphicx}
\usepackage[svgnames]{xcolor}
\usepackage[normalem]{ulem}

\numberwithin{equation}{section}

\usepackage{wrapfig}
\newcommand{\pv}{\mbox{p}_V}
\newcommand{\xt}{\widetilde{x}}
\newcommand{\yt}{\widetilde{y}}

\usepackage[all]{xy}

\usepackage{enumitem}
\usepackage{xcolor}
\definecolor{darkcandyapplered}{rgb}{0.64, 0.0, 0.0}
\definecolor{joliblanc}{RGB}{252, 251, 249}
\definecolor{bluegray}{rgb}{0.4, 0.6, 0.8}
\definecolor{ceruleanblue}{rgb}{0.16, 0.32, 0.75}
\definecolor{greydark}{RGB}{75, 71, 59}
\definecolor{oldlavender}{RGB}{75, 71, 59}

\colorlet{shadecolor}{white}

\usepackage{hyperref}
\hypersetup{colorlinks=true,linkcolor=darkcandyapplered,citecolor=oldlavender}

\usepackage{caption}

\def\defn#1{{\bf\itshape #1}}

\usepackage{thmtools}

\declaretheoremstyle[%
  spaceabove=6pt,%
  spacebelow=6pt,%
  headfont=\normalfont\itshape,%
  postheadspace=1em,%
  qed=\qedsymbol,%
  headpunct={}
]{mystyle} 
\declaretheorem[name={Proof ---},style=mystyle,unnumbered,
]{Proof}

\renewcommand{\qedsymbol}{\color{oldlavender}{$\blacksquare$}}

\usepackage{amsthm}

\newtheoremstyle{mytheoremstyle} 
    {\topsep}                    
    {\topsep}                    
    {\itshape}                   
    {}                           
    {\scshape}                   
    {.}                          
    {.5em}                       
    {}  

\theoremstyle{mytheoremstyle}

\newtheorem{theorem}{Theorem}[section]
\newtheorem{proposition}[theorem]{Proposition}
\newtheorem{cor}[theorem]{Corollary}
\newtheorem{lemma}[theorem]{Lemma}

\newtheoremstyle{mydefinitionstyle} 
    {\topsep}                    
    {\topsep}                    
    {}                   
    {}                           
    {\itshape}                   
    {.}                          
    {.5em}                       
    {}  

\theoremstyle{mydefinitionstyle}

\newtheorem{definition}[theorem]{Definition}
\newtheorem*{remark}{Remark}
\newtheorem{example}[theorem]{Example}

\usepackage{geometry}
\usepackage{titlesec}
\usepackage{secdot}
\titleformat*{\section}{\centering\bfseries}
\titleformat*{\subsection}{\itshape}

\usepackage{enumitem}

\usepackage{bbm}

\newcommand{\T}{\mbox{T}}
\newcommand{\SO}{SO} 
\newcommand{\SE}{SE} 

\renewcommand{\O}{O}

\newcommand{\R}{\mathbb{R}}
\newcommand{\Z}{\mathbb{Z}}

\newcommand{\C}{\mathbb{C}}

\newcommand{\ad}{\mbox{ad}}
\newcommand{\Ad}{\mbox{Ad}}

\newcommand{\J}{\Phi}

\newcommand\restr[2]{{
  \left.\kern-\nulldelimiterspace 
  #1 
  \vphantom{\big|} 
  \right|_{#2} 
  }}
 \newcommand{\IB}{\mathbb{I}_{\mathcal{B}}}
\newcommand{\IF}{\mathbb{I}_{\mathcal{F}}}

\newcommand{\TG}{T}

\renewcommand{\NG}{N}

\renewcommand{\a}{a}
\newcommand{\x}{x}
\newcommand{\y}{y}
\renewcommand{\b}{b}

\newcommand{\0}{0}

\newcommand{\cat}{\mbox{Cat}}
  
\newcommand{\partialt}{\left. \frac{\mbox{d}}{\mbox{dt}}\right| _{t=0}}

\renewcommand{\d}{\mathrm{d}}
\newcommand{\ddt}{\frac{\d}{\d t}}

\usepackage{authblk}
\author{\small\textsc{Marine Fontaine} and \textsc{James Montaldi}}
\date{}

\title{Persistence of stationary motion under explicit symmetry breaking perturbation}

\begin{document}
\maketitle

\begin{abstract}
Explicit symmetry breaking occurs when a dynamical system having a certain symmetry group is perturbed to a system which has strictly less symmetry. We give a geometric approach to study this phenomenon in the setting of hamiltonian systems.   We provide a method for determining the equilibria and relative equilibria that persist after a  symmetry breaking perturbation.  In particular a lower bound for the number of each is found, in terms of the equivariant Lyusternik-Schnirelmann category of the group orbit.
\vspace{6pt}

\noindent \small\textit{Keywords:} Symmetry breaking, hamiltonian systems, Lie group actions
\end{abstract}


\section{Introduction} 
When we talk about symmetries, we either refer to the symmetry of a physical law (dynamical equations) or the symmetry of a physical state (solution of these equations). The \defn{symmetry} or \defn{symmetry group }of a physical law (or a physical state) is defined to be the group of transformations which leave these equations (or this solution) invariant. \defn{Explicit symmetry breaking} is defined as a process of perturbing symmetric dynamical equations so that the resulting equations have a lower symmetry group. When the system is not hamiltonian interesting results are obtained showing for example that periodic solutions of an unperturbed dynamical system  can become heteroclinic cycles under a perturbation that breaks the symmetry \cite{MR2092068,MR1850318,MR1187861}. 

In this paper we focus on dynamical systems which are hamiltonian.  We address the question whether equilibria or relative equilibria persist under a symmetry breaking perturbation. 

Some aspects of explicit symmetry breaking phenomena for hamiltonian systems have been studied by several authors \cite{MR879691,Grabsi,MR2684603} to cite just a few, but except for \cite{Grabsi} their results do not overlap with ours (for \cite{Grabsi} see further below).  In the case when none of the symmetries are broken the question of persistence of relative equilibria is raised in \cite{MR1438262} for compact symmetry groups and further developed in \cite{MR1660367, Wulff}. The persistence in those papers was for perturbing the momentum value rather than the Hamiltonian function, but the arguments also apply if the Hamiltonian is perturbed in a manner that preserves the symmetry. 
 
In applications explicit symmetry breaking phenomena appear in various ways. As explained for example in \cite{Brading:706822} terms can be introduced artificially in the equations of motion in order to match with theoretical or experimental observations. Besides quantization processes might also be a cause for the appearance of such terms which are the so-called quantum anomalies \cite{MR3202835,MR3791131}. In this case the terms are not artificially introduced but they appear after a renormalization procedure.

The problem is as follows. Phase spaces of hamiltonian systems are symplectic manifolds and the symmetries of such systems are encoded into Lie group actions on those manifolds. A symplectic manifold is a smooth manifold $M$ equipped with a non-degenerate closed two-form $\omega$. A (proper) smooth action of a Lie group $G$ on $M$ is \defn{symplectic} if it preserves $\omega$. An important class of symplectic group actions on symplectic manifolds are the \defn{Hamiltonian} actions, which are those actions to which there is  associated a Noether conserved quantity expressed in term of a momentum map $\J_G:M\to\mathfrak{g}^*$, where $\mathfrak{g}^*$ is the dual of the Lie algebra of $G$. This notion generalizes the notion of angular momentum in classical mechanics, when the phase space is (a product of copies of) $T^*\R^3$, acted on by the group of rotations $SO(3)$. By a \defn{hamiltonian (proper) $G$-manifold} we mean a quadruple $(M,\omega,G,\J_G)$ as described above, where $\J_G:M\to\mathfrak{g}^*$ is equivariant with respect to the coadjoint action $\Ad^*:(g,\mu)\in G\times\mathfrak{g}^*\mapsto\Ad^*_{g^{-1}}\mu\in\mathfrak{g}^*$ of $G$ on the dual Lie algebra.

The dynamics is governed by a Hamiltonian $h$ which is a $G$-invariant smooth real-valued function defined on $M$. The ring of such functions is denoted by $C^{\infty}(M)^G$. The non-degeneracy of $\omega$ implies that, associated to any Hamiltonian $h\in C^{\infty}(M)^G$, there is a unique vector field $X_h$ defined by $\iota_{X_h}\omega=-\d h$. Since the action of $G$ on $M$ is symplectic and $h$ is $G$-invariant, the integral curve $\varphi_t(m)$ of $X_h$ starting at $m\in M$ satisfies $\varphi_t(g\cdot m)=g\cdot \varphi_t(m)$ for all $g\in G$. The resulting Hamiltonian equations
\begin{equation}\label{original}
		\ddt\varphi_t(m)=X_h(\varphi_t(m))
\end{equation}
are thus $G$-equivariant and we say that $G$ is the \defn{symmetry group} of \eqref{original}. We study the effect of a small hamiltonian perturbation of these equations, which is invariant with respect to a subgroup of $G$.
\begin{definition} \label{def:H-perturbation}
\normalfont Let $h\in C^{\infty}(M)^G$ and $H\subset G$ be a closed subgroup. An \defn{$H$-pertubation} of $h$ is a family of functions $h_{\lambda}\in C^{\infty}(M)^H$, with $\lambda$ in a neighbourhood of 0 in $\R$ such that the map $(m,\lambda)\in M\times\R\mapsto h_{\lambda}(m)\in \R$ is smooth and $h_0=h$.
\end{definition}

 We focus on specific solutions of \eqref{original}, namely equilibria (fixed points under the dynamics) and relative equilibria (group orbits fixed under the dynamics). Under a specific non-degeneracy condition on a (relative) equilibrium of the unperturbed Hamiltonian $h$ there is a chance that this (relative) equilibrium persists under an $H$-perturbation. 

Section \ref{s: Persistence of equilibria} is devoted to the question of persistence of equilibria. The required non-degeneracy condition on an equilibrium $m\in M$ of $h$ is a particular case of the Morse-Bott condition when the critical manifold of $h$ is the group orbit $G\cdot m$ (cf. Definition \ref{morse bott non degenerate}). We show that at least a certain number of $H$-orbits of equilibria persist under a small $H$-perturbation in a tubular neighbourhood of $G\cdot m$ (cf. Theorem \ref{crit} and Corollary \ref{persistence of equilibria}). This number is the positive integer $\cat_H\left(G/G_m\right)$, which is the $H$-equivariant Lyusternik-Schnirelmann category of the group orbit. At the end of the section we present applications of this result including the problem of an ellipse-shaped planar rigid body moving in a planar irrotational, incompressible fluid with zero vorticity and zero circulation around the body.

Extending Theorem \ref{crit} and Corollary \ref{persistence of equilibria} to the case of relative equilibria is more subtle because we must take into account the conservation of momentum and the corresponding non-degeneracy condition takes into account the ambient symplectic structure. This question is treated in Section \ref{s: symmetry breaking for RE}. Whereas equilibria are critical points of the Hamiltonian function $h$, relative equilibria are critical points of the restriction of this same function to a level set of the momentum map, the problem being that as the group changes, so do these level sets. Let $m\in M$ be one of those critical points. The element $\xi\in\mathfrak{g}$ playing the role of a Lagrange multiplier is called the velocity of $m$, which is in general not unique when the action is not free. For that reason we refer to a relative equilibrium as a pair $(m,\xi)\in M\times \mathfrak{g}$.  
We denote the underlying Lagrange function associated to $\xi$ by $h^\xi = h-\phi_G^\xi$, where $\phi_G^{\xi}(m):=\langle \Phi_G(m),\xi\rangle$ is a $G_{\xi}$-invariant function on $M$. It is called the \emph{augmented Hamiltonian}.

A standard definition says that a relative equilibrium $(m,\xi)$ of $h$ is non-degenerate if the Hessian of $h^{\xi}$ at $m$ is a non-singular quadratic form when restricted to some symplectic subspace $\NG_1\subset T_mM$, called the symplectic slice at $m$. If the perturbations $h_{\lambda}$ are invariant with respect to the full symmetry group $G$, this notion of non-degeneracy is enough to guarantee the persistence of a relative equilibrium. This is no longer the case if $h_{\lambda}$ has a smaller symmetry group than does $h$ and we require a stronger non-degeneracy condition on the relative equilibrium (Definition \ref{alpha-nondegeneracy}). In \cite{Grabsi} a step in that direction is taken, when the symmetry group is a torus that breaks into a subtorus. In addition, the group actions in consideration are assumed to be free. We extend their result to non-free actions and non-abelian symmetry groups. A necessary condition for a relative equilibrium of $h$ to persist under an $H$-perturbation is that the velocity $\xi$ belongs to $\mathfrak{h}$, the Lie algebra of $H$. If the non-degeneracy condition on $(m,\xi)\in M\times\mathfrak{h}$ holds, and modulo some technicalities, the least number of $H_{\mu}$-orbits of relative equilibria with velocity close to $\xi$, which persist under a small $H$-perturbation in some neighbourhood of $G_{\mu}\cdot m$ in $\J_H^{-1}(\alpha)$, is the positive integer $\cat_{H_{\mu}}\left(G_{\mu}/G_m\right).$
This is the content of Theorem \ref{persistence of RE} and Corollary \ref{persistence of RE cor}. We illustrate this result for the spherical pendulum on $S^3$, as a perturbation of the geodesic flow; an example of symmetry breaking from $SO(4)$ to $SO(3)$.

\paragraph{Acknowledgments.} We would like to thank Luis Garc\'ia-Naranjo for suggesting the example of the $2$D-rigid body submerged in a fluid, also the editor and referees for suggestions improving the exposition. This work forms part of the first author's Ph.D.\ thesis \cite{moi} from the University of Manchester. It was partially funded by the project ``symplectic techniques in differential geometry'' within the Excellence of Science program of the F.R.S.-FNRS and FWO.

\section{Symmetry breaking for equilibria}\label{s: Persistence of equilibria}
The aim of this section is to give a lower bound for the number of $H$-orbits of equilibria that persist under a small $H$-perturbation of some $G$-invariant Hamiltonian $h$.  Since equilibria of the hamiltonian vector field are the same as critical points of the Hamiltonian function, we state the main theorem in terms of critical points of smooth functions.  The persistence result will require a non-degeneracy condition which we now recall.

\begin{definition}\label{morse bott non degenerate}
	\normalfont A \defn{$G$-nondegenerate critical point} of an invariant function  $h\in C^{\infty}(M)^G$ is a point $m\in M$ such that
	\begin{enumerate}[label=(\roman*)]
			\item $\d h(m)=0$,\label{morse bott i}
			\item if $\mathcal{N}$ is any subspace of $T_mM$ complementary to $\mathfrak{g}\cdot m$, the restriction $D_{\mathcal{N}}^2h(m)$ of the Hessian $D^2h(m)$ to $\mathcal{N}\times \mathcal{N}$ is non-singular. In other words, the Hessian is non-singular in the directions normal to the group orbit. \label{non degenerate hessian}
	\end{enumerate}
\end{definition}

\begin{remark}
\begin{enumerate}[label=(\roman*)]
\item We also say an equilibrium of an invariant Hamiltonian is $G$-nondegenerate under the same conditions. 
\item If $m\in M$ is a $G$-nondegenerate equilibrium of $h$ then so is any $p\in G\cdot m$, by $G$-invariance. For this reason, the tangent space $T_p\left(G\cdot m\right)$ is contained in $\ker\left(D^2h(p)\right)$ for any $p\in G\cdot m$. Definition \ref{morse bott non degenerate} is a particular case of Morse-Bott non-degeneracy when $G\cdot m$ is the critical manifold of $h$ (cf. \cite{Bott}). Note that Condition \ref{non degenerate hessian} implies that the critical manifold $G\cdot m$ is isolated in the sense that there exists a tubular neighbourhood of $G\cdot m$ that does not contain any other critical points of $h$. 
\item  This definition is also valid in infinite dimensional Hilbert spaces, provided the non-degeneracy is interpreted as saying that the linear map $\mathcal{N}\to\mathcal{N}^*$ given by $\mathbf{v}\to D_{\mathcal{N}}^2h(m)(\mathbf{v})$ is \emph{invertible} (has bounded inverse).

\end{enumerate}
\end{remark}

\subsection{Persistence of critical points}
We say that a closed subgroup $H\subset G$ is \defn{co-compact} (in $G$) if the left multiplication of $H$ on $G$ is co-compact, i.e.\ the orbit space $H\setminus G$ under this action is compact.

\begin{theorem}\label{crit}
	Let $G$ be a Lie group acting properly on a manifold $M$ and let $H\subset G$ be a closed co-compact subgroup. Assume that $h_{\lambda}\in C^{\infty}(M)^H$ is an $H$-pertubation of some $h\in C^{\infty}(M)^G$ in the sense of Definition\,\ref{def:H-perturbation}, and that $m\in M$ is a $G$-nondegenerate equilibrium of $h$. 
		
	Then there is a $G$-invariant neighbourhood $U\subset M$ of $m$ such that, if $\lambda$ is sufficiently small, there exists a function $f_{\lambda}\in C^{\infty}(G/G_m)^H$ whose critical points are in one-to-one correspondence with those of $h_{\lambda}$ in $U$.
\end{theorem}

\def\Ne{\mathcal{N}_\varepsilon}
\begin{Proof} 
 Let $m\in M$ be a $G$-nondegenerate equilibrium of $h$ whose stabilizer is denoted by $K:=G_m$.  Let $\mathcal{N}\subset T_mM$ be a $K$-invariant vector subspace complementary to $\mathfrak{g}\cdot m$ in $T_mM$.  Recall that (see e.g., \cite{MR2021152}) the natural action of the direct product $G\times K$ on $G\times \mathcal{N}$ is given by 
\begin{equation*}
	(h,k)\cdot (g,\nu)=(hgk^{-1},k\cdot \nu).
\end{equation*} 
The Palais Tube Theorem then states that there is a $G$-invariant neighbourhood $U\subset M$ of $m$ and a $K$-invariant neighbourhood $\Ne$ of $0$ in $\mathcal{N}$ such that the associated bundle $G\times_K\Ne\subset G\times_K\mathcal{N}$ is a local model for $U$ and the only critical points of $h$ in $U$ are on $G\cdot m$. In that model the point $m$ reads $[(e,0)]$ and the $H$-pertubation is identified with a family also denoted $h_{\lambda}:G\times_K\Ne\to \R$. Let $\rho:G\times \Ne\to G\times_K\Ne$ be the orbit map for the $K$-action on $G\times\Ne$. We define the lift of $h_{\lambda}$ by 
		$$\widetilde{h}_{\lambda} := h_{\lambda}\circ\rho:G\times \Ne\to \R.$$   
		The critical points of $\widetilde{h}_\lambda$ are then the inverse image under $\rho$ of those of $h_{\lambda}$.
We may thus work with $\widetilde{h}_{\lambda}$ instead of $h_{\lambda}$. 

By assumption the lift $\widetilde{h}$ is $G\times K$-invariant and $\widetilde{h}_{\lambda}$ is $H\times K$-invariant. Since $(e,0)\in G\times \Ne$ is a $G$-nondegenerate critical point of $\widetilde{h}$,
\begin{equation}\label{hknot}
				\d \widetilde{h}(e,0)=0\quad\mbox{and}\quad D^2_{\mathcal{N}}\widetilde{h}(e,0)\quad\mbox{is non-singular}.
\end{equation}
In particular the map 
\begin{equation*}
	\d_{\mathcal{N}}\widetilde{h}:G\times \Ne\to \mathcal{N}^*
\end{equation*} 
satisfies $\d_{\mathcal{N}}\widetilde{h}(e,0)=0$ and its derivative with respect to the $\mathcal{N}$-variables, evaluated at $(e,0)$, is invertible.  

We now wish to apply the implicit function theorem globally and uniformly along the orbit $G\times\{0\}$.  The implicit function theorem implies the existence, for each $g\in G$, of a neighbourhood $V_g\times W_g$ of $(0,g)$ in $\R\times G$ such that, for each $(\lambda,g')\in V_g\times W_g$, there is a unique $\phi^g_{\lambda}(g')\in \Ne$ satisfying 
\begin{equation}\label{dn}
		\d_{\mathcal{N}}\widetilde{h}_{\lambda}(g',\phi^g_{\lambda}(g'))=0.
\end{equation}
for all $g'\in W_g$.
By $H\times K$-invariance of $\widetilde{h}_{\lambda}$ and hence of $\d_{\mathcal{N}}\widetilde{h}_{\lambda}$, we can choose $W_g$ to be $H\times K$-invariant. Note that when we refer to $H$ or $K$ individually we are thinking of them as subgroups of $H\times K$.
This procedure defines an $H$-invariant smooth function 
\begin{eqnarray*}
\phi^g:V_g\times W_g & \longrightarrow & \Ne\\
(\lambda,g') & \longmapsto & \phi^g_{\lambda}(g').
\end{eqnarray*} 
Since the $W_g$ are $H$-invariant open subsets of $G$, and $H\setminus G$ is compact, it follows that we can extract a finite subcover of $\{W_g\}$, call this $\{W_{g_1},\dots,W_{g_n}\}$. Now let $V_0=V_{g_1}\cap\dots\cap V_{g_n}$. By restricting the above maps to $\lambda\in V_0$, we obtain a globally (in $G$) defined map 
$$\phi: V_0\times G \longrightarrow \Ne$$
satisfying $\d_{\mathcal{N}}\widetilde{h}_\lambda(g,\phi_\lambda(g))=0$, for all $(\lambda,g)\in G\times V_0$ (this gives the global and uniform application of the implicit function theorem we required). Moreover, $\phi$ is $H$-invariant and $K$-equivariant:
\begin{equation}\label{eq:phi_lambda invariance}
\phi_\lambda(hgk^{-1}) = k\cdot\phi_\lambda(g).
\end{equation}

Now define a family $f_\lambda$ of functions on $G$ by
\begin{equation}\label{f_lambda}
f_\lambda(g) = \widetilde{h}_\lambda(g,\phi_\lambda(g)).
\end{equation}
We claim that $f_\lambda$ is $H\times K$-invariant and has a critical point at $g$ if and only if $\widetilde{h}_\lambda$ has a critical point at $(g,\phi_\lambda(g))\in G\times\Ne$. The $H\times K$-invariance implies that $f_\lambda$ passes down to a smooth $H$-invariant function on $G/K$, as required.

To check the critical point property, note that, with $(g,w)=(g,\phi_\lambda(g))$,
$$\d f_\lambda(g) = \d_G \widetilde{h}_\lambda(g,w) + \d_{\mathcal{N}}\widetilde{h}_\lambda(g,w)\d\phi_\lambda(g).$$
However, for $\lambda\in V_0$, $w=\phi_\lambda(g)$ if and only if \eqref{dn} holds, with $\phi$ in place of $\phi^g$.  Thus $f_\lambda$ has a critical point at $g$ if and only if $\widetilde{h}_\lambda$ has a critical point at $(g,\phi_\lambda(g))$.

The invariance properties of $f_\lambda$ follows from those of $\phi_\lambda$ given in \eqref{eq:phi_lambda invariance} above.
\end{Proof}

\subsection{Persistence of critical points and equilibria in hamiltonian systems}
Before progressing to give a lower bound for the number of critical points of the perturbation, we recall an important concept in the calculus of variations. In their original paper~\cite{LS}, Lyusternik and Schnirelmann introduce a numerical homotopy invariant of a topological space $M$ that they denote $\cat(M)$. They define it to be the least number of open subsets of $M$, whose inclusion is nullhomotopic, that are required to cover $M$. They show that if $M$ is a closed (i.e. compact without boundary) smooth manifold, then any smooth function $f$ on $M$ has at least $\cat(M)$ critical points. 
The equivariant analogue $\cat_G(M)$ when $G$ is a compact Lie group is obtained in \cite{fadell,marzantowicz} and the extension to proper Lie group actions in \cite{ayala}. The \defn{equivariant Lyusternik-Schnirelmann category}, denoted $\cat_G(M)$, is the least number of $G$-invariant open subsets of $M$, contractible by mean of a $G$-homotopy onto a $G$-orbit, that are required to cover $M$. The extension of the result of \cite{LS} to proper $G$-manifolds (possibly infinite dimensional) requires a certain compactness condition, called the \defn{orbitwise Palais-Smale condition (OPS)} \cite{ayala}. In our applications $M$ is finite dimensional and the orbit space $M/G$ is compact, and in this case the OPS condition is automatic. 

\begin{theorem}[Equivariant Lyusternik-Schnirelmann Theorem \cite{ayala,bartsch}]\label{gcrit}
If a proper $G$-manifold $M$ and a function $f\in C^{\infty}(M)^G$ satisfy condition (OPS), then $f$ has at least $\cat_G(M)$ group orbits of critical points.
\end{theorem}

 As a corollary of Theorem \ref{crit} we obtain:

\begin{cor}
	Under the assumptions of Theorem \ref{crit} the number of $H$-orbits of critical points of $h$ that persist near $G\cdot m$ under a small $H$-perturbation is bounded below by $\cat_H(G/G_m)$.
\end{cor}

\begin{Proof}
	If $\lambda$ is sufficiently small, Theorem \ref{crit} implies that the $H$-orbits of critical points of $h_{\lambda}$ in some neighbourhood of $G\cdot m$ are in one-to-one correspondence with those of a function $f_{\lambda}\in C^{\infty}(G/K)^H$ where $K:=G_m$. Since $H\setminus G$ is compact it follows that so too is $H\setminus G/K$. We may therefore apply Theorem \ref{gcrit} and conclude that the number of $H$-orbits of critical points of $h_{\lambda}$ is at least $\cat_H(G/K)$ near $G\cdot m$.
\end{Proof}

\begin{remark}\label{rmk:morse}
If one knows \emph{a priori} that the critical points of $f_\lambda$ are all $H$-nondegenerate, then a (usually better) lower bound can be given using equivariant Morse theory, see for example \cite{MR702806,MR739783}. 
\end{remark}

The corollary can be reformulated in the hamiltonian setting.

\begin{cor}[Persistence of Equilibria]\label{persistence of equilibria}
Suppose $h\in C^{\infty}(M)^G$ is a Hamiltonian defined on a symplectic proper $G$-manifold $(M,\omega)$ and has a $G$-nondegenerate equilibrium $m$. Then the number of $H$-orbits of equilibria that persist near $G\cdot m$ under a small $H$-perturbation is bounded below by $\cat_H(G/G_m)$.
\end{cor}

\begin{example}\label{first example}
\normalfont 
Think of the cylinder $M=S^1\times \R$ as embedded in $\R^3$ with coordinates $(\theta,z)$ and endow it with the standard symplectic form $\omega=\d \theta\wedge \d z$. The Lie group $G=\O(2)$ acts on $M$ by $R_{\varphi}\cdot (\theta,z)=(\theta+\varphi, z),$ if $R_{\varphi}\in \O(2)$ is a rotation of angle $\varphi$; and by $r_{\alpha}\cdot (\theta,z)=(2\alpha-\theta, z),$ if $r_{\alpha}\in \O(2)$ is a reflection about the line forming an angle $\alpha$ with the $x$-axis in $\R^3$. The action of $G$ on $M$ is hamiltonian with momentum map $\J_G:(\theta,z)\in M\mapsto z\in \R$. Consider the $1$-parameter family $h_{\lambda}:S^1\times \R\to \R$ defined by $$h_{\lambda}(\theta,z)=z^2+\lambda\cos(n\theta).$$ Then $h=h_0$ is $G$-invariant and $m=(0,0)$ is a $G$-nondegenerate equilibrium of $h$ whose stabilizer is $G_m=\langle r_0\rangle$. The perturbation $h_{\lambda}$ is invariant under  $H=D_n$, where $D_n$ is the dihedral group of order $2n$. In fact, the full symmetry group is $D_n\times \Z_2$ since $\Z_2$ acts on the $z$-component by changing its sign. However such an action is not symplectic. Since this discrete part does not contribute in the further application, we do not take it into account. The perturbed Hamiltonian $h_{\lambda}$ has $2n$ critical points whose coordinates are $(\frac{\pi}{n}k,0)$ for $k=0,\ldots,2n-1$, which form a regular $2n$-gone as shown in Figure \ref{fig:cylinder_image} for the case $n=3$. Since $G/G_m=\O(2)/\langle r_0\rangle$ is topologically a circle, we find $\cat_H(G/G_m)=2$ (cf. \cite{marzantowicz} Corollary $1.17$). There are thus two $H$-orbits of equilibria of $h$ which persist, each of them being a regular $n$-gone (cf. Figure \ref{fig:orbitofequilibria_image}).

\begin{figure}[!ht]
	\centering
	\begin{minipage}[t]{7cm}
		\centering
		\includegraphics [width =2.5 cm]{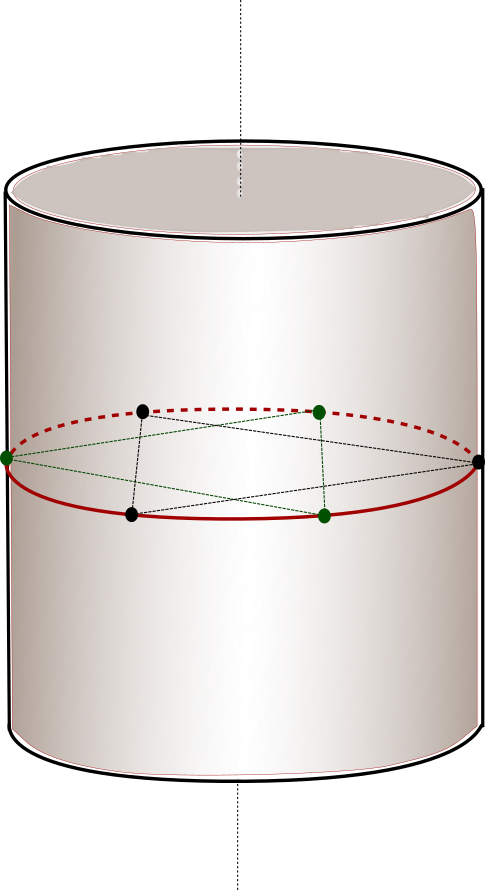}
		\caption{\small When $n=3$, $h$ has a $G$-orbit (red circle) consisting of $G$-nondegenerate equilibria, on which the six equilibria of $h_{\lambda}$ lie.}
		\label{fig:cylinder_image}
	\end{minipage}\hfill
	\begin{minipage}[t]{7cm}
		\centering
		\includegraphics [width =5.5 cm]{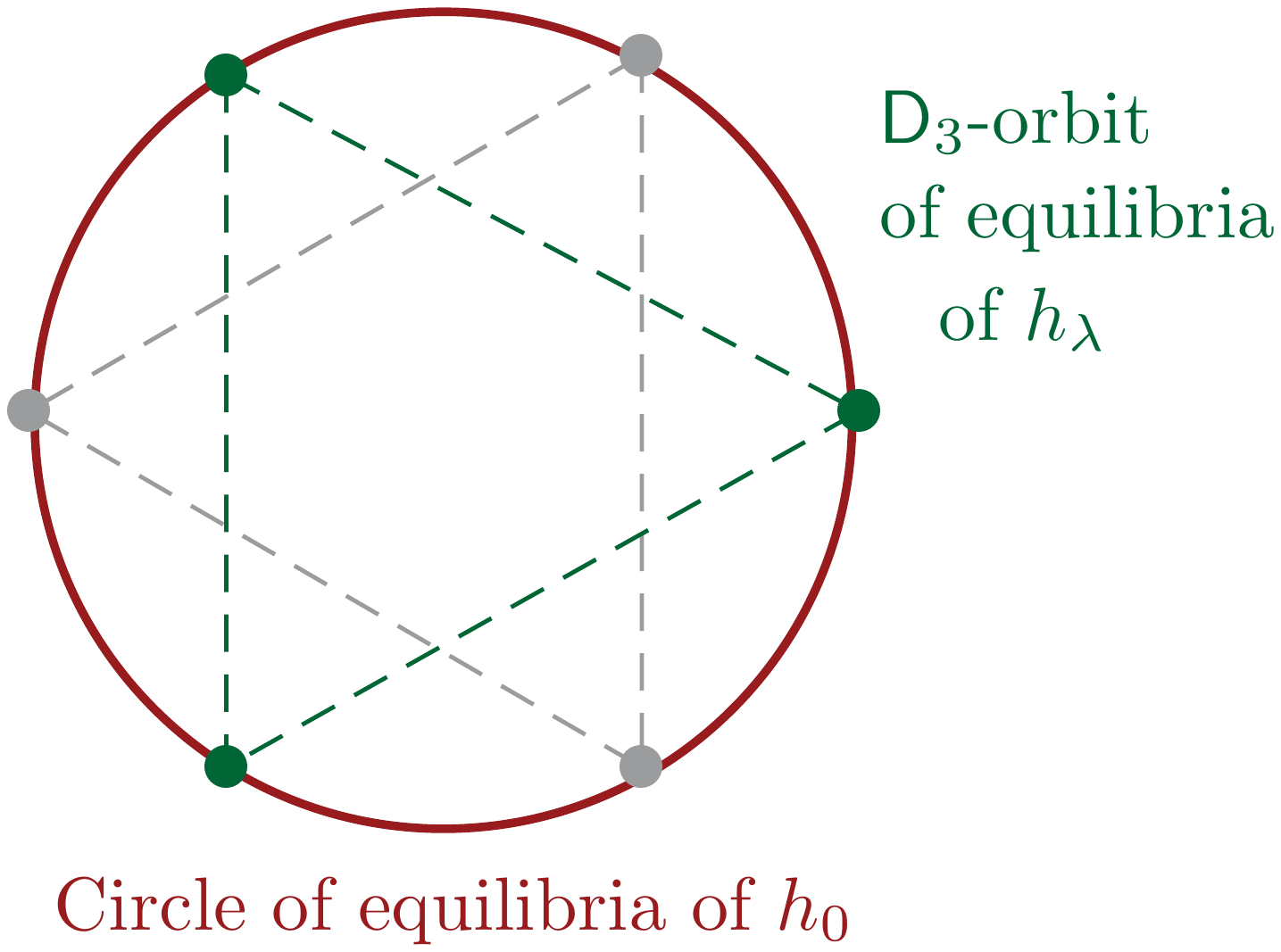}
		\caption{\small At the level of coordinate $z=0$, the six equilibria of $h_{\lambda}$ form two different $D_3$-orbits. One orbit is stable and one is unstable.}
		\label{fig:orbitofequilibria_image}
	\end{minipage}
\end{figure}
\end{example}

\subsection{Dynamics of a 2D rigid body in a potential flow}\label{ex: body-fluid}
We apply the result of Corollary \ref{persistence of equilibria} to the problem of a planar rigid body $\mathcal{B}$ of mass $m$ moving in a planar irrotational, incompressible fluid with zero vorticity and zero circulation around the body. The motion is governed by Kirchhoff equations~\cite{Kirchhoff}. Classical treatments of the problem can be found in \cite{Lamb} and \cite{Milne-Thomson}. The configuration space of the body-fluid system is a submanifold $Q$ of the product $SE(2)\times\mbox{Emb}_{vol}\left(\mathcal{F}_0,\R^2\right)$, where $SE(2)$ is the special Euclidean group describing the motion of the body, and $\mbox{Emb}_{vol}\left(\mathcal{F}_0,\R^2\right)$ is the space of volume-preserving embeddings of the fluid reference space $\mathcal{F}_0$ in $\R^2$. The symmetry group of this system is the direct product of $SE(2)$ (group of uniform body-fluid translations and rotations) and the particle relabeling symmetry group (volume-preserving diffeomorphisms of $\mathcal{F}_0$). Since these actions commute, the system can be reduced by the process of symplectic reduction by stages \cite{MR2337886}.

The Hamiltonian of the system is invariant under the particle relabeling symmetry group. Geometrically, eliminating the fluid variables amounts to carry out a symplectic reduction by this group. The particle relabeling symmetry group acts on $T^*Q$ in a hamiltonian fashion. The associated momentum map has two components corresponding to the vorticity and the circulation. The reduction at zero momentum corresponds to a fluid with zero circulation and zero vorticity. In this case, the symplectic reduced space is identified with $T^*\SE(2)$, endowed with the canonical symplectic form and the $\SE(2)$-invariant reduced Hamiltonian is the sum of the kinetic energy of the body-fluid system by the addition of the so-called ``added masses'', and the kinetic energy of the body. Those added masses depend only on the body's shape and not on the mass distribution. The reader is refered to \cite{Kanso} and \cite{fluid} for details. Since $\SE(2)$ acts symplectically on $T^*\SE(2)$, the dynamics can be reduced a second time using Poisson reduction and thereby the reduced motion is governed by the Kirchhoff equations that are the Lie-Poisson equations on the dual Lie algebra $\mathfrak{se}(2)^*$.

For the sake of simplicity we will assume that the body $\mathcal{B}$ is shaped as an ellipse with semi-axes of length $A>B>0$. We will use the formulae and follow the notations of \cite{Luis}. At the center of mass of $\mathcal{B}$ we attach a frame $\left\lbrace E_1,E_2\right\rbrace$ that is aligned with the symmetry axes of the body. Its position is related at any time to a fixed space frame $\left\lbrace e_1,e_2\right\rbrace$ by an element of $\SE(2)$. An element of the Lie algebra $\xi\in\mathfrak{se}(2)$ is identified with a vector 
\begin{equation}\label{xi}
	(\dot{\theta},v_1,v_2)\in\R^3
\end{equation} 
where $\dot{\theta}\in \R$ is the angular velocity of $\mathcal{B}$ and $(v_1,v_2)^T\in\R^2$ is the linear velocity of its center of mass, expressed in the body's frame. In this setting the body has kinetic energy
\begin{equation}
	T_{\mathcal{B}}=\frac{1}{2}\xi\cdot \IB\xi
\end{equation}
with $\IB:=\mbox{diag}(I_{\mathcal{B}},m,m)$, where $I_{\mathcal{B}}$ is the moment of inertia of the body about its center of mass. The kinetic energy of the fluid is given by 
\begin{equation}
	T_{\mathcal{F}}=\frac{1}{2}\xi\cdot \IF\xi
\end{equation}
where $\IF=\frac{\rho\pi}{4} \mbox{diag}((A^2-B^2)^2,B^2,A^2)$ is the tensor of added masses, and $\rho$ is the fluid density. In the absence of external forces, the Lagrangian of the body-fluid system $\mathcal{L}:T\SE(2)\to \R$ is given by $\mathcal{L}=T_{\mathcal{B}}+T_{\mathcal{F}}$. It defines a Riemannian metric on $\SE(2)$ with respect to which the motion of the body $\mathcal{B}$ is geodesic. Since $\mathcal{L}$ does not depend on the group variables, it is $\SE(2)$-invariant and can thus be reduced to the function $\ell:\mathfrak{se}(2)\to\R$ given by
\begin{equation}
	\ell(\xi)=\frac{1}{2}\xi\cdot(\IB+\IF)\xi
\end{equation}
with $\xi$ as in \eqref{xi}. An element $\nu$ of the dual Lie algebra $\mathfrak{se}(2)^*$ is identified with a one by three matrix $(x,\alpha_1,\alpha_2)$. The dual pairing $\langle\cdot ,\cdot\rangle$ between $\mathfrak{se}(2)^*$ and $\mathfrak{se}(2)$ is thus given by
\begin{equation}
	\langle \nu,\xi\rangle:=(x,\alpha_1,\alpha_2)(\dot{\theta},v_1,v_2)^T=x\dot{\theta}+\alpha_1 v_1+\alpha_2 v_2.
\end{equation}
We perform the Legendre transform $\mathbb{F}L:\xi\in\mathfrak{se}(2)\mapsto ((\IB+\IF)\xi)^T\in\mathfrak{se}(2)^*$ to obtain the reduced Hamiltonian $h:\mathfrak{se}(2)^*\to \R$ defined by $$h(\nu)=\frac{1}{2}\nu\cdot(\IB+\IF)^{-1}\nu^T.$$ The Lie-Poisson equations on $\mathfrak{se}(2)^*$ that describe the motion of the body-fluid system are 
\begin{equation}\label{Lie poisson}
	\dot{\nu}=\mbox{ad}^*_{\frac{\delta h}{\delta \nu}}\nu.
\end{equation}
where $\ad_{\xi}^*\nu$ is identified with $(\alpha_1v_2-\alpha_2v_1,\dot{\theta}\alpha_2,-\dot{\theta}\alpha_1)$. This problem turns out to exhibit symmetry breaking phenomena from different points of view:
\begin{enumerate}[label=(\roman*)]
	\item One point of view consists in looking at the body $\mathcal{B}$ without the fluid ($\rho=0$). Adding the fluid amounts to seeing the fluid density $\rho$ as a ``parameter''. The $\O(2)$-symmetry of the kinetic reduced Hamiltonian breaks into a $D_2$-symmetry, where $D_2$ is the symmetry group of an ellipse.
	\item On the other hand we can consider the original system as being a circular planar rigid body ($A=B$) in a fluid and the symmetry can be broken by deforming the body into an elliptical shaped body. This case exhibits the same pattern of symmetry breaking from $\O(2)$ to the subgroup $D_2$.
\end{enumerate} 
These two approaches are the same from a group theoretical point of view. Contrary to Example \ref{first example}, the Hamiltonian in consideration will not be perturbed by adding some potential energy. In this case, there is no potential energy involved, only the metric is perturbed giving rise to a modified kinetic energy. Let us now discuss the two cases mentioned above.
\begin{enumerate}[label=(\roman*)]
	\item The unperturbed system on the Poisson reduced space $\mathfrak{se}(2)^*$ is governed by the Hamiltonian 
\begin{equation}
	h(\nu)=\frac{1}{2}\nu\cdot\mathbb{I}\nu=\frac{1}{2}\left(\frac{x^2}{I_{\mathcal{B}}}+\frac{\alpha_1^2+\alpha_2^2}{m}\right)
\end{equation} 
where $\nu:=(x,\alpha_1,\alpha_2)$ and $\mathbb{I}:=\IB^{-1}$. The Hamiltonian is invariant with respect to the group $G=O(2)$.  In particular, for each $c\in\R$, the level sets $h(\nu)=c$ describe spheroids in $\R^3$. Note that the full symmetry group is in fact $O(2)\times \Z_2$ since $\Z_2$ acts on the $x$-component by swapping the sign. However this discrete part does not contribute to our analysis.

Adding a fluid to the system amounts to look at the variation of the dimensionless parameter
$$\lambda=d\rho\quad\mbox{where}\quad d:=\frac{A^2-B^2}{m}>0\quad\mbox{is fixed}.$$ This gives rise to the perturbed Hamiltonian $h_{\lambda}(\nu)=\frac{1}{2}\nu\cdot \mathbb{I}_{\lambda}\nu$ with 
\begin{equation}
\mathbb{I}_{\lambda}=\mbox{diag}\left(\frac{1}{I_{\mathcal{B}}+\lambda c_1},\frac{1}{m+\lambda c_2},\frac{1}{m+\lambda c_3}\right).
\end{equation}
where $c_1=\frac{m^2d\pi}{4}, c_2=\frac{\pi(A^2-m d)}{4d}$ and $c_3=\frac{\pi(B^2+m d)}{4d}$ are fixed constants encoding the datas of the system.
The perturbed Hamiltonian reads
\begin{equation}\label{h ex}
h_{\lambda}(\nu)=\frac{1}{2}\left(\frac{x^2}{I_{\mathcal{B}}+\lambda c_1}+\frac{\alpha_1^2}{m+\lambda c_2}+\frac{\alpha_2^2}{m+\lambda c_3}\right)
\end{equation}
and has symmetry $H=D_2$, the dihedral group of order four: recall that the group $D_2$ is isomorphic to $\Z_2\times\Z_2$, and acts here by changing the signs of $\alpha_1$ and $\alpha_2$. 

This perturbation coincides with $h$ when $\lambda=0$ and the function $(\lambda,\nu)\mapsto h_{\lambda}(\nu)$ is smooth. Therefore, $h_{\lambda}$ is an $H$-pertubation of $h$. The symmetry is broken because the fluid influences the motion of the body if it is elliptical. If the body is circular ($A=B$), or if it moves in the vacuum, its center of mass would move at constant velocity and it would rotate at constant angular speed.

\item We carry out another kind of perturbation: rather than perturbing the rigid body motion by adding a fluid to the system, we start with a circular planar rigid body ($A=B$) in a fluid and break the symmetry by changing the body shape into an ellipse. The unperturbed Hamiltonian is given by 
\begin{equation}
	h(\nu)=\frac{1}{2}\nu\cdot\mathbb{I}\nu=\frac{1}{2}\left(\frac{x^2}{I_{\mathcal{B}}}+\frac{\alpha_1^2+\alpha_2^2}{m+d_2}\right)
\end{equation} 
where $d_2=\frac{\rho\pi B^2}{4}$, $\nu:=(x,\alpha_1,\alpha_2)$, $\mathbb{I}:=(\IB+\IF)^{-1}$ and $A=B$ in the definition of $\IF$. The Hamiltonian is invariant with respect to $G=\O(2)$. For each $c\in\R$, the level sets $h(\nu)=c$ also describe spheroids in $\R^3$.

We perturb the body shape by setting $\lambda=\frac{A^2-B^2}{B^2}$ where $B>0$ is fixed and $A\geq B>0$ varies. This gives rise to the perturbed Hamiltonian $h_{\lambda}(\nu)=\nu\cdot \mathbb{I}_{\lambda}\nu$ with 
\begin{equation}
\mathbb{I}_{\lambda}=\mbox{diag}\left(\frac{1}{I_{\mathcal{B}}+\lambda^2 d_1},\frac{1}{m+d_2},\frac{1}{m+(\lambda+1) d_2}\right)
\end{equation}
where $d_1=\frac{\rho\pi B^4}{4}$.
The perturbed Hamiltonian is thus given by
\begin{equation}
h_{\lambda}(\nu)=\frac{1}{2}\left(\frac{x^2}{I_{\mathcal{B}}+\lambda^2 d_1}+\frac{\alpha_1^2}{m+d_2}+\frac{\alpha_2^2}{m+(\lambda+1) d_2}\right)
\end{equation}
and is again symmetric with respect to the action of $H=D_2$. In this case, if there was no fluid ($\rho=d_2=0)$, no symmetries would have been broken.
\end{enumerate}

\begin{figure}[!ht]
	\centering	
		\begin{minipage}[t]{9cm}
		\centering
		\includegraphics [width =6 cm]{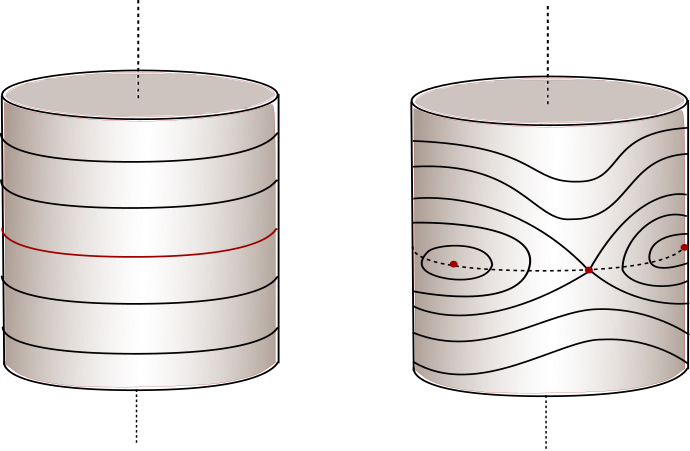}		
		\caption{\small The flow lines are given by intersecting the level sets of $h_{\lambda}$ (the ellipsoids) and the coadjoint orbits. On the left hand side, we see the flow lines of $h$ on a coadjoint orbit. On the right hand side, the flow has been perturbed.}\label{fig:flow lines}
\end{minipage}	
\end{figure}

Since the reduced motion is governed by the Lie-Poisson equations \eqref{Lie poisson}, it is constrained to the coadjoint orbits of $\SE(2)$. As shown in \cite{MR} (Chapter $14.6$), almost all of them are cylinders (the singular orbits consist of points on the vertical dashed line in Figure \ref{fig:flow lines}). In both cases, the level sets of $h_{\lambda}$ are ellipsoids and those of $h=h_0$ are spheroids. Their intersections with a coadjoint orbit are shown in Figure \ref{fig:flow lines}. In particular, the circle of equilibria of $h$ (in red in Figure \ref{fig:flow lines}) breaks into four fixed points of $h_{\lambda}$, two of which are connected by four heteroclinic cycles.
\vspace{0.5cm}

Let us go back to the first case we discussed above with $h_{\lambda}$ as in \eqref{h ex}. We will apply Corollary \ref{persistence of equilibria} to predict the existence of the four fixed points that persist (cf. Figure \ref{fig:flow lines}). The Fr\'echet derivative of $h_{\lambda}$ is 
\begin{equation}
\frac{\delta h_{\lambda}}{\delta \nu}=\left(\frac{x}{I_{\mathcal{B}}+\lambda c_1},\frac{\alpha_1}{m+\lambda c_2},\frac{\alpha_2}{m+\lambda c_3}\right).
\end{equation} 
Therefore, the Lie-Poisson equations \eqref{Lie poisson} reduces to
\begin{equation}\label{lp}
\left\lbrace
\begin{array}{ccc}
\dot{x} & = & \frac{\lambda(c_2-c_3)}{(m+\lambda c_3)(m+\lambda c_2)}\alpha_1\alpha_2\\
\dot{\alpha_1} & = & \frac{x\alpha_2}{m+\lambda c_1}\\
\dot{\alpha_2} & = & -\frac{x\alpha_1}{m+\lambda c_1}
\end{array}\right.
\end{equation}
Setting $\lambda=0$ in \eqref{lp}, we see that the fixed points of $h=h_0$ are either of the form $(0,\alpha_1,\alpha_2)$ with $(\alpha_1,\alpha_2)\in(\R^2)^*$, or of the form $(x,0,0)$ which correspond to points on the singular coadjoint orbit.

Let $\mu:=(0,\alpha_1,\alpha_2)$ with $\alpha_1^2+\alpha_2^2=1$ be a fixed point of the unperturbed hamiltonian $h$. The isotropy subgroup of $\mu$ is $G_{\mu}=\langle r_{\vartheta}\rangle$ where $r_{\vartheta}$ is a reflection in the plane. The quotient $G/G_{\mu}=\O(2)/\langle r_{\vartheta}\rangle$ is topologically a circle yielding $\cat_{D_2}(S^1)=2.$ The four fixed points appearing in Figure \ref{fig:flow lines} are the two $H$-orbits that persist.

\section{Symplectic preliminaries} \label{s: the symplectic slice theorem}
For our principal result in Section\,\ref{s: symmetry breaking for RE} below, the symplectic geometry is crucial, and we rely strongly on the symplectic local model for a hamiltonian proper $G$-manifold near a group orbit.

The Symplectic Tube Theorem is used to study the local dynamics and the local geometry of a hamiltonian proper $G$-manifold $(M,\omega,G,\J_G)$. It states essentially that every $m\in M$ admits a $G$-invariant neighbourhood, which is $G$-equivariantly symplectomorphic to a neighbourhood of the zero section of a symplectic associated bundle. This construction provides tractable semi-global coordinates for $M$ near $G$-orbits. Those coordinates are sometimes referred as \defn{slice coordinates}. This theorem appears in \cite{Guillemin, Marle85} for symplectic Lie group actions with equivariant momentum maps.  Its extension to general symplectic Lie group actions can be found in \cite{MR2021152,MR1486529}. See also \cite{schmah,Perlmutter} for the particular case where $M$ is a cotangent bundle.

We briefly recall the construction underlying the Symplectic Tube Theorem. The reader is referred to \cite{MR2021152} or \cite{MR3242761} for details. Let $m\in M$ with momentum $\mu=\J_G(m)$. Denote by $G_m$ and $G_{\mu}$ the stabilizers of $m$ and $\mu$ respectively and by $\mathfrak{g}_m$ and $\mathfrak{g}_{\mu}$ their respective Lie algebras. The stabilizer $G_m$ is compact by properness of the action of $G$ on $M$. We can thus split $\mathfrak{g}_{\mu}$  and $\mathfrak{g}$ into a direct sum of $G_m$-invariant subspaces 
\begin{equation*}
	\mathfrak{g}_{\mu}=\mathfrak{g}_m\oplus \mathfrak{m}\quad \mbox{ and }\quad\mathfrak{g}=\mathfrak{g}_m\oplus \mathfrak{m}\oplus \mathfrak{n}.
\end{equation*} 
We denote by $\mathfrak{g}\cdot m$ the tangent space at $m$ of $G\cdot m$. Elements of $\mathfrak{g}\cdot m$ are vectors of the form $x_M(m):=\partialt \mbox{exp}(tx)\cdot m$, where $x\in\mathfrak{g}$ and $\mbox{exp}:\mathfrak{g}\to G$ is the group exponential. The tangent space $T_mM$ can be decomposed into a direct sum of four $G_m$-invariant subspaces
\begin{equation}\label{Witt Artin G}
	T_mM=\TG_0\oplus \TG_1\oplus \NG_0\oplus \NG_1,
\end{equation}
known as the Witt-Artin decomposition, defined as follows:
\begin{enumerate}[label=(\roman*)]
		\item $\TG_0:=\ker\left(D\J_G(m)\right)\cap \mathfrak{g}\cdot m=\mathfrak{g}_{\mu}\cdot m$.		
		\item $\TG_1:=\mathfrak{n}\cdot m$ which is a symplectic vector subspace of $(T_mM,\omega(m))$.
		\item $\NG_1$ is a choice of $G_m$-invariant complement to $\TG_0$ in $\ker\left(D\J_G(m)\right)$. It is a symplectic subspace of $(T_mM,\omega(m))$  and is called the \defn{symplectic slice}. The linear action of $G_m$ on $\NG_1$ is hamiltonian with momentum map $\J_{\NG_1}:\NG_1\to\mathfrak{g}_m^*$ given by $\langle\J_{\NG_1}(\nu),x\rangle=\frac{1}{2}\omega(x_{N_1}(\nu),\nu)$ for every $\nu\in\NG_1$ and $x\in\mathfrak{g}_m$.
		\item $\NG_0$ is a $G_m$-invariant Lagrangian complement to $\TG_0$ in the symplectic orthogonal $(\TG_1\oplus \NG_1)^{\omega(m)}$. There is an isomorphism $f:\NG_0\to \mathfrak{m}^*$ given by $\langle f(w),y\rangle=\omega(m)\left(y_M(m),w\right)$ for every $w\in \NG_0$ and $y\in \mathfrak{m}$.
	\end{enumerate}
Since $\NG_1$ is a $G_m$-invariant subspace, there is a well-defined action of $G_m$ on the product $G\times \mathfrak{m}^*\times \NG_1$ given by 
\begin{equation}\label{twisted action}
	k\cdot (g,\rho,\nu)=(gk^{-1},\Ad^*_{k^{-1}}\rho,k\cdot\nu).
\end{equation}
This action is free and proper by freeness and properness of the action on the $G$-factor. The orbit space $Y$ is thus a smooth manifold whose points are equivalence classes of the form $[(g,\rho,\nu)]$. The group $G$ acts smoothly and properly on $Y$, by left multiplication on the $G$-factor. Let $\mathfrak{m}^*_{\varepsilon}\subset \mathfrak{m}^*$ and $(\NG_1)_{\varepsilon}\subset \NG_1$ be $G_m$-invariant neighbourhoods of zero in $\mathfrak{m}^*$ and $\NG_1$, respectively. Then
\begin{equation}\label{tube}
	Y_{\varepsilon}:=G\times_{G_m}(\mathfrak{m}^*_{\varepsilon}\times (\NG_1)_{\varepsilon})
\end{equation}
is a neighbourhood of the zero section in $Y$. It comes with a symplectic form $\omega_{Y_{\varepsilon}}$ if it is chosen small enough \cite{MR2021152} (Proposition $7.2.2$).

\begin{theorem}[Symplectic Tube Theorem]\label{symplectic slice thm}
	Let $(M,\omega,G,\J_G)$ be a hamiltonian proper $G$-manifold. Let $m\in M$ with momentum $\mu=\J_G(m)$ and let $(Y_{\varepsilon},\omega_{Y_{\varepsilon}})$ as in \eqref{tube}. Then there is a $G$-invariant neighbourhood $U\subset M$ of $m$ and a $G$-equivariant symplectomorphism $\varphi: (Y_{\varepsilon},\omega_{Y_{\varepsilon}})\to (U,\restr{\omega}{U})$ such that $\varphi\left([e,0,0]\right)=m$.
\end{theorem}

\noindent We call the triplet $(\varphi,Y_{\varepsilon},U)$ a \defn{symplectic $G$-tube} at $m$ and we also say that $(Y_{\varepsilon},\omega_{Y_{\varepsilon}})$ is a \defn{symplectic local model} for $(U,\restr{\omega}{U})$. Moreover, the momentum map $\J_G:M\to\mathfrak{g}^*$ can be expressed in terms of the slice coordinates:

\begin{theorem}[Marle-Guillemin-Sternberg Normal Form]\label{MGS}
	Let $(M,\omega,G,\J_G)$ be a hamiltonian proper $G$-manifold and let $(\varphi,Y_{\varepsilon},U)$ be a symplectic $G$-tube at $m\in M$. Then the $G$-action on $Y_{\varepsilon}$ is hamiltonian with associated momentum map $\widetilde{\Phi}_G:Y_{\varepsilon}\to \mathfrak{g}^*$ defined by
\begin{equation}\label{Jyr}
	\widetilde{\Phi}_G([g,\rho,\nu])=\Ad^*_{g^{-1}}(\J_G(m)+\rho+\J_{N_1}(\nu)).
\end{equation}
If $G$ is connected, \eqref{Jyr} coincides with $\restr{\J_G}{U}$ when pulled back along $\varphi^{-1}$.
\end{theorem}

\section{Symmetry breaking for relative equilibria}\label{s: symmetry breaking for RE}

In this section, we extend Theorem \ref{crit} and Corollary \ref{persistence of equilibria} to the case of relative equilibria which is more subtle for two reasons: firstly we must take into account the conservation of momentum, and secondly for a non-zero velocity the so-called augmented Hamiltonian no longer has symmetry $G$.

We start by briefly recalling some standard facts about relative equilibria (see \cite{lecturesonmechanics} for details). Given a hamiltonian proper $G$-manifold $(M,\omega,G,\J_G)$, a \defn{relative equilibrium} of a Hamiltonian $h\in C^{\infty}(M)^G$ is a pair $(m,\xi)\in M\times\mathfrak{g}$ such that $X_h(m)=\xi_M(m)$. Equivalently, if $(m,\xi)$ is a relative equilibrium of $h$, then $m$ is a critical point of the \defn{augmented Hamiltonian} $$h^{\xi}:=h-\phi_G^{\xi}\in C^{\infty}(M)^{G_{\xi}}$$ where $\phi_G^{\xi}(m):=\langle \J_G(m),\xi\rangle$, which is a $G_{\xi}$-invariant function which depends linearly on $\xi$. A standard fact about relative equilibria is that the velocity $\xi$ and the momentum $\mu=\J_G(m)$ commute i.e. $\xi\in\mathfrak{g}_{\mu}$.
Note that, if the isotropy group $G_m$ is non trivial and $(m,\xi)$ is a relative equilibrium of $h$, then $(m,\xi+\eta)$ is also a relative equilibrium of $h$, for any $\eta\in\mathfrak{g}_m$. Moreover if $(m,\xi)$ is a relative equilibrium of $h$ then so is $(g\cdot m,\Ad_g\xi)$ for every $g\in G$.  In general a relative equilibrium is said to be non-degenerate if the Hessian $D^2h^{\xi}(m)$ is a non-singular quadratic form, when restricted to the symplectic slice $\NG_1$ at $m$ relative to the $G$-action. However, this definition of non-degeneracy is not enough to guarantee that a relative equilibrium of some $h\in C^{\infty}(M)^G$ persists under an $H$-perturbation. For that reason, we need a stronger version of non-degeneracy.

\subsection{Induced momentum map}
Let $H$ be a closed subgroup of $G$. The dual of the inclusion of Lie algebras $i_{\mathfrak{h}}:\mathfrak{h}\hookrightarrow \mathfrak{g}$ is the projection $\mathfrak{i}_{\mathfrak{h}}^*:\mathfrak{g}^*\to \mathfrak{h}^*$ and is given by $i_{\mathfrak{h}}^*(\mu)=\restr{\mu}{\mathfrak{h}}$, which is the restriction of the linear form $\mu$ to the Lie subalgebra $\mathfrak{h}$. The action of $H$ on $M$ is still both symplectic and Hamiltonian. A momentum map for this action is given by $\J_H=i_{\mathfrak{h}}^*\circ\J_G:M\to\mathfrak{h}^*$ and is called the \defn{induced momentum map} for the $H$-action. 

\begin{proposition}[\cite{slicesubgroup}]\label{result A}
Consider the decomposition of $T_mM$ as in \eqref{Witt Artin G}, and define the subspace 
$$\mathcal{M}:=\left\lbrace z_M(m)+w\in \TG_1\oplus\NG_0\mid -\ad_z^*\mu+f(w)\in\mathfrak{h}^{\circ}\right\rbrace$$ 
where $f$ denotes the isomorphism between $\NG_0$ and $\mathfrak{m}^*$, and $\mathfrak{h}^{\circ}$ is the annihilator of $\mathfrak{h}$ in $\mathfrak{g}^*$. Then $\ker\left(D \J_H(m)\right)= \ker\left(D \J_G(m)\right)\oplus \mathcal{M}.$
\end{proposition}
\begin{Proof}
It is clear from the definitions that there is an inclusion of subspaces
\begin{equation}\label{kerg in kerk}
		\ker\left(D\J_G(m)\right)\subset \ker\left(D\J_H(m)\right).
\end{equation} Let $\left(\varphi, G\times_{G_m} \left(\mathfrak{m}^*_{\varepsilon}\times (\NG_1)_{\varepsilon}\right),U\right)$ be a symplectic $G$-tube at $m$ as in Theorem \ref{symplectic slice thm}. Linearising $\varphi^{-1}$ at $m$ yields a linear symplectomorphism
\begin{equation*}\label{Tphi}
		T_m\varphi^{-1}: \TG_0\oplus \TG_1\oplus \NG_0\oplus \NG_1 \to T_{\varphi^{-1}(m)}\left( G\times_{G_m} \left(\mathfrak{m}^*\times \NG_1\right)\right).
\end{equation*}
For $x+y\in\mathfrak{g}_{m}\oplus \mathfrak{m}$ and $z\in \mathfrak{n}$ we have
\begin{equation*}
		T_m\varphi^{-1}\cdot ((x+y)_M(m)+z_M(m)+w+\nu)=T_{(e,0,0)}\rho\cdot (x+y+z,f(w),\nu)
\end{equation*}
where $\rho:G\times \mathfrak{m}^*\times \NG_1\to G\times_{G_m}\left(\mathfrak{m}^*\times \NG_1\right)$ is the orbit map. By definition, the subspace $\ker\left(D\J_H(m)\right)$ consists of the elements $$((x+y)_M(m)+z_M(m)+w+\nu)\in \TG_0\oplus \TG_1\oplus \NG_0\oplus \NG_1$$ satisfying $D(\restr{\J_H}{U}\circ \varphi\circ \rho)(e,0,0)\cdot (x+y+z,f(w),\nu)=0.$ Equivalently 
\begin{eqnarray*}
		0 &=& \partialt \restr{\J_H}{U}\circ \varphi\left([(\exp({t(x+y+z)}),tf(w),t\nu)]\right)\\
		&=& \partialt i^*_{\mathfrak{h}}\left(\Ad^*_{\exp(-t(x+y+z))}\left(\mu+tf(w)+\J_{\NG_1}(t\nu)\right)+C\right),\qquad C\in \R\\
		&=&i^*_{\mathfrak{h}}\left(-\ad^*_z\mu+f(w)\right)\\
\end{eqnarray*}
where the normal form for the momentum map is given by Theorem \ref{MGS}. As required $-\ad^*_z\mu+f(w)\in\mathfrak{h}^{\circ}$ since the kernel of $i^*_{\mathfrak{h}}$ is equal to $\mathfrak{h}^{\circ}$. Note that we do not need the assumption of Theorem \ref{MGS} that $G$ is connected because the statement only depends on the differential.
\end{Proof}
\subsection{Non-degeneracy condition and regularity condition}
We now state a stronger version of non-degeneracy of a relative equilibrium.
\begin{definition}\label{alpha-nondegeneracy}
	\normalfont Let $(M,\omega,G,\J_G)$ be a hamiltonian proper $G$-manifold with $H\subset G$ be a closed subgroup, and $\J_H:M\to \mathfrak{h}^*$ be the induced momentum map. Setting $\alpha:=\J_H(m)$, a relative equilibrium $(m,\xi)\in M\times \mathfrak{g}$ of $h\in C^{\infty}(M)^G$ is said to be \defn{$\alpha$-nondenegerate} if $D^2h^{\xi}(m)$ is a non-singular quadratic form on $\NG_1\oplus \mathcal{M}$ with $\mathcal{M}$ as in Proposition \ref{result A}.
\end{definition}
Definition \ref{alpha-nondegeneracy} only depends on $\alpha$ and not on the underlying Witt-Artin decomposition of $T_mM$. If $G$ is non-abelian, the space $\mathcal{M}$ might have an non-trivial intersection with $\mathfrak{g}\cdot m$. This intersection is the subspace $\mathfrak{q}\cdot m\subset \mathfrak{g}\cdot m$ where $\mathfrak{q}$ is an $H_m$-invariant complement to $\mathfrak{g}_{\mu}$ in the ``symplectic orthogonal'' 
$$\mathfrak{h}^{\perp_{\mu}}:=\left\lbrace x\in\mathfrak{g}\mid x_M(m)\in(\mathfrak{h}\cdot m)^{\omega(m)}\right\rbrace.
$$
The non-singularity of $D^2h^{\xi}(m)$ along $\mathfrak{g}\cdot m$ depends only on that of $D^2 \phi^{\xi}_G(m)$ which has symmetry group $G_{\xi}$. The condition is a consequence of the following lemma which is proved in Section \ref{section hessian}.

\begin{lemma}\label{hessian degeneracy}
		Let $(M,\omega,G,\J_G)$ be a hamiltonian proper $G$-manifold. Let $m\in M$ with momentum $\mu=\J_G(m)$ and an element $\xi\in \mathfrak{g}_{\mu}$. If $\mathfrak{g}$ is semi-simple then the Hessian $D^2\phi_G^{\xi}(m)$ restricted to $\mathfrak{g}\cdot m$ is singular precisely along $(\mathfrak{g}_{\xi}+\mathfrak{g}_{\mu})\cdot m$.  
\end{lemma}
Therefore if an equilibrium $(m,\xi)\in M\times\mathfrak{g}$ of some $h\in C^{\infty}(M)^G$ with momentum $\mu=\J_G(m)$ is $\alpha$-nondegenerate in the sense of Definition \ref{alpha-nondegeneracy}, then $\mathfrak{g}_{\xi}$ has trivial intersection with $\mathfrak{q}$. In Theorem \ref{persistence of RE} we show that a number of orbits of relative equilibria of $h$ persist under $H$-perturbations. Such relative equilibria must have their velocity $\xi$ in $\mathfrak{h}_{\mu}$. We assume an additional regularity assumption 
\begin{equation}\label{assumption}\tag{R}
	\mathfrak{g}_{\mu}\subset\mathfrak{g}_{\xi}
\end{equation} 
This says essentially that $\mu$ needs to be more regular than $\xi$ (cf. \cite{moi} Definition $6.2.2$ for more details).  In particular if condition (R) is satisfied then the null-space of the Hessian referred to in Lemma \ref{hessian degeneracy} above is equal to $\mathfrak{g}_\xi\cdot m$.

\begin{example}\label{ex:so(3) in so(4)}
	\normalfont In this example we show when condition \eqref{assumption} holds for $\mathfrak{g}=\mathfrak{so}(4)$ and subalgebras $\mathfrak{h}$ isomorphic to $\mathfrak{so}(3)$. The Lie algebra $\mathfrak{g}$ is identified with the set of pairs $(x,a)\in\R^3\times\R^3$ with Lie bracket
\begin{equation}\label{lie bracket so(4)}
	[(\x,\a),(\y,\b)]=(\x\times\y+\a\times \b,\x\times\b+\a\times \y).
\end{equation}
The dual Lie algebra $\mathfrak{g}^*$ consists of pairs $(\chi,\rho)\in\R^3\times\R^3$ which satisfy
\begin{equation*}
	\langle (\chi,\rho),(\x,\a)\rangle=\chi\cdot\x+\rho\cdot \a.
\end{equation*}
The linearized coadjoint action of $\mathfrak{g}$ on $\mathfrak{g}^*$ is given by 
\begin{equation}\label{coadjoint so(4)}
	\ad^*_{(\x,\a)}(\chi,\rho)=(\chi\times\x+\rho\times\a,\chi\times\a+\rho\times\x).
\end{equation}

\noindent\textit{Lie subalgebras isomorphic to $\mathfrak{so}(3)$.}  Elements of $\mathfrak{so}(3)$ are identified with vectors $\x\in\R^3$. We consider two inequivalent Lie subalgebras of $\mathfrak{h}\subset\mathfrak{g}$ isomorphic to $\mathfrak{so}(3)$, namely

\begin{enumerate}[label=(\roman*)]
		\item The Lie algebra of rotations in $\R^3$ denoted
		$\mathfrak{so}(3)_r=\left\lbrace (\x,\0)\in \R^6\mid \x\in \R^3\right\rbrace$ with Lie bracket $$[(\x,\0),(\y,\0)]=(\x\times\y,\0).$$
		\item The diagonal elements denoted $\mathfrak{so}(3)_d=\left\lbrace \left(\frac{\x}{2},\frac{\x}{2}\right)\in \R^6\mid \x\in \R^3\right\rbrace$ with Lie bracket $$[(\frac{\x}{2},\frac{\x}{2}),(\frac{\y}{2},\frac{\y}{2})]=(\frac{\x\times\y}{2},\frac{\x\times\y}{2}).$$
\end{enumerate}

\noindent\textit{Regularity condition.} Given a fixed momentum $\mu:=(\chi,\rho)\in \mathfrak{g}^*$, the stabilizer Lie subalgebra is
\begin{equation*}
	\mathfrak{g}_{\mu}=\left\lbrace (x,a)\in \mathfrak{g}\mid \chi\times x+\rho\times a=0\quad\mbox{and}\quad\chi\times a+\rho\times x=0\right\rbrace
\end{equation*}
by \eqref{coadjoint so(4)}. We show below whether condition \eqref{assumption} is satisfied for our two different choices of Lie subalgebra.

\begin{enumerate}[label=(\roman*)]
\item Let $\mathfrak{h}=\mathfrak{so}(3)_r$ with inclusion map $i_{\mathfrak{h}}:\x\in\mathfrak{h}\mapsto (\x,\0)\in\mathfrak{g}.$ To compute the dual of this inclusion $i_{\mathfrak{h}}^*:\mathfrak{g}^*\to \mathfrak{h}^*$, we take $(\chi,\rho)\in \mathfrak{g}^*$ and $\x\in\mathfrak{h}$ and we compute
\begin{equation*}
	\langle i_{\mathfrak{h}}^*(\chi,\rho),\x\rangle=\langle (\chi,\rho),i_{\mathfrak{h}}(\x)\rangle=\langle (\chi,\rho),(\x,\0)\rangle=\chi\cdot \x.
\end{equation*} 
Then $$i_{\mathfrak{h}}^*((\chi,\rho))=\chi\in \mathfrak{h}^*.$$

The symplectic orthogonal is $\mathfrak{h}^{\perp_{\mu}}=\left\lbrace (x,a)\in \mathfrak{g}\mid \chi\times x+\rho\times a=0\right\rbrace.$ Since the velocity $\xi\in\mathfrak{h}$ must commute with $\mu$, it has to belong to the subspace $\mathfrak{h}_{\mu}=\mathfrak{g}_{\mu}\cap\mathfrak{h}$. Using equation \eqref{coadjoint so(4)},
\begin{equation*}
		\mathfrak{h}_{\mu}=\left\lbrace (x,0)\in\mathfrak{so}(3)_r\mid \chi\times x=0\mbox{ and }\rho\times x=0\right\rbrace.
\end{equation*}
There are three cases to consider:
\begin{enumerate}[label=(\alph*)]
\item If $\chi=\rho=0$ then $\mathfrak{g}_{\mu}=\mathfrak{g}$ and $\mathfrak{h}_{\mu}=\mathfrak{h}$. We choose $\xi=(y,0)\in\mathfrak{h}$ where $y\in\R^3$ is arbitrary. Using \eqref{lie bracket so(4)} we get
$$\mathfrak{g}_{\xi}=\left\lbrace (\lambda_1 y,\lambda_2 y)\in \mathfrak{g}\mid \lambda_1,\lambda_2\in\R\right\rbrace$$ and clearly \eqref{assumption} does not hold.

\item If $\chi$ and $\rho$ are not collinear, $\mathfrak{h}_{\mu}=\left\lbrace (0,0)\right\rbrace$. In this case, the only available velocity is $\xi=0$ and thus $\mathfrak{g}_{\xi}=\mathfrak{g}$. In particular \eqref{assumption} holds.
\item If $\mu=(\chi,\rho)$ is such that $\chi=s\rho$ for some $s\in \R$, we choose $\xi\neq0$ of the form $$\xi:=(\lambda\chi,0)\in\mathfrak{h}_{\mu}\quad\mbox{for some}\quad\lambda\in\R$$ 
and thus $\mathfrak{g}_{\xi}=\left\lbrace (x,a)\in \mathfrak{g}\mid x\times\chi=0\mbox{ and }a\times \chi=0\right\rbrace.$ Note that in particular, $\mathfrak{g}_{\xi}\subset\mathfrak{g}_{\mu}$. To see whether $\mathfrak{g}_{\mu}\subset \mathfrak{g}_{\xi}$, pick an element $(x,a)\in \mathfrak{g}_{\mu}$. By definition, it satisfies
\begin{equation}\label{eq}
	x\times\chi=\rho\times a\quad\mbox{and}\quad \chi\times a=x\times\rho.
\end{equation}
Using \eqref{eq} and the fact that $\chi=s\rho$ we get,
\begin{equation*}
	x\times\chi=s(x\times \rho)=s(\chi\times a)=s^2(\rho\times a)=s^2(x\times\chi).
\end{equation*}
Similarly
\begin{equation*}
	a\times\chi=s(a\times\rho)=s(\chi\times x)=s^2(\rho\times x)=s^2(a\times\chi).
\end{equation*}
Therefore, $(x,a)\in\mathfrak{g}_{\xi}$ as long as $s^2\neq 1$; that is,  \eqref{assumption} holds for such $\xi$ as long as $\mu\neq(\chi,\pm\chi)$ (see Figure \ref{R condition figure}). 
\end{enumerate}

\begin{figure}[!ht]
	\centering	
		\begin{minipage}[t]{8cm}
		\centering
		\includegraphics [width =4cm]{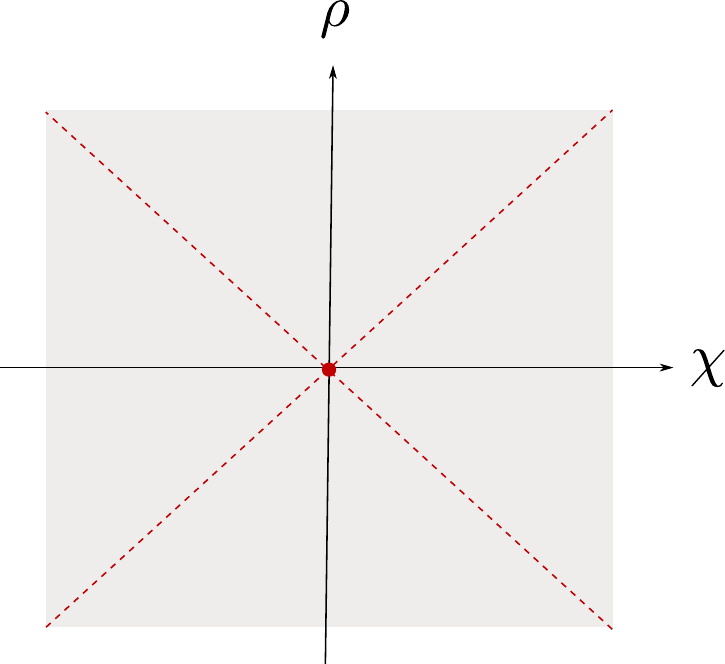}		
		\caption{\small Condition \eqref{assumption} holds as long as $\mu$ is away form the red dashed lines which represent subspaces of codimension three in $\R^6$.}\label{R condition figure}
\end{minipage}	
\end{figure}
\item Let $\mathfrak{h}=\mathfrak{so}(3)_d$ with inclusion map $$i_{\mathfrak{h}}:\x\in\mathfrak{h}\mapsto \left(\frac{\x}{2},\frac{\x}{2}\right)\in\mathfrak{g}.$$ To compute the dual of this inclusion $i_{\mathfrak{h}}^*:\mathfrak{g}^*\to \mathfrak{h}^*$, we take $(\chi,\rho)\in \mathfrak{g}^*$ and $\x\in\mathfrak{h}$ and we compute
\begin{equation*}
	\langle i_{\mathfrak{h}}^*(\chi,\rho),\x\rangle=\langle (\chi,\rho),i_{\mathfrak{h}}(\x)\rangle=\langle (\chi,\rho),\left(\frac{\x}{2},\frac{\x}{2}\right)\rangle=\frac{\chi+\rho}{2}\cdot \x,
\end{equation*} 
Then
\begin{equation*}
		i_{\mathfrak{h}}^*((\chi,\rho))=\frac{\chi+\rho}{2}\in \mathfrak{h}^*.
\end{equation*} 
Set $\mu:=(\chi,\rho)\in \mathfrak{g}^*$ and $\alpha:=i_{\mathfrak{h}}^*(\mu)=\frac{\chi+\rho}{2}\in \mathfrak{h}^*$. Using Equation \eqref{coadjoint so(4)} we get $$\mathfrak{h}_{\mu}=\left\lbrace \left(\frac{x}{2},\frac{x}{2}\right)\in\mathfrak{so}(3)_d\mid \alpha\times x=0\right\rbrace.$$
We thus choose a velocity of the form $$\xi:=(\lambda \alpha,\lambda \alpha)\in\mathfrak{h}_{\mu}$$ for some $\lambda\in\R$. By \eqref{lie bracket so(4)} the stabilizer Lie algebra of $\xi$ is
\begin{equation}\label{gxi}
		\mathfrak{g}_{\xi}=\left\lbrace (x,a)\in \mathfrak{g}\mid x\times\alpha+a\times \alpha=0\right\rbrace.
\end{equation}
In particular, $\mathfrak{g}_{\mu}\subset \mathfrak{g}_{\xi}$ and \eqref{assumption} is satisfied for any choice of $\mu$.
\end{enumerate}
\end{example}

\subsection{Persistence of relative equilibria}\label{sub: persistence of RE}
We are now ready to state the corresponding version of Theorem \ref{crit} for relative equilibria. The proof follows the same steps as Theorem \ref{crit}. For that reason some details have been skipped.

\begin{theorem}\label{persistence of RE}
	Let $(M,\omega,G,\J_G)$ be a hamiltonian proper $G$-manifold. Assume that $h\in C^{\infty}(M)^G$ has a relative equilibrium $(m,\xi)\in M\times \mathfrak{h}$ with momentum $\mu=\J_G(m)$. Let $\alpha\in\mathfrak{h}^*$ be the restriction of $\mu$ to $\mathfrak{h}$. We assume that
	\begin{enumerate}[label=(\roman*)]
		\item $\alpha$ is a regular value of $\Phi_H$,
		\item $(m,\xi)$ is $\alpha$-nondegenerate and assumption \eqref{assumption} is satisfied,
		\item $H_{\mu}\subset G_{\mu}$ is co-compact.
	\end{enumerate}

Then there is a $G_{\mu}$-invariant neighbourhood $U\subset \J_H^{-1}(\alpha)$ of $m$ and a neighbourhood $V\subset\R\times\mathfrak{h}$ of $(0,\xi)$ such that, for each $(\lambda,\eta)\in V$, there is a function $f^{\eta}_{\lambda}\in C^{\infty}(G_{\mu}/G_m)^{H_{\mu}}$, depending smoothly on $(\lambda,\eta)\in V$,  whose critical points are in one-to-one correspondence with those of $h^{\eta}_{\lambda}$ in $U$. 
\end{theorem}

\begin{Proof} The proof is similar to that of Theorem \ref{crit}, except that $h^\xi$ is not $G$-invariant, but only $G_\xi$ (or thanks to condition \eqref{assumption} $G_\mu$-invariant), and the extra ingredient is to control the difference between $G$ and $G_\mu$-nondegeneracy, and between the level sets of $\Phi_G$ and $\Phi_H$. 

	Let $(m,\xi)\in M\times \mathfrak{h}$ be an $\alpha$-nondegenerate relative equilibrium $h$, where $\alpha$ is the restriction of the momentum $\mu=\J_G(m)$ to $\mathfrak{h}$. In particular $\xi\in\mathfrak{h}_{\mu}$. 
	Set $K:=G_m$ and fix a Witt-Artin decomposition \eqref{Witt Artin G} relative to the $G$-action. Note that $N_1\oplus T_1$ is a symplectic slice at $m$ relative to the $G_{\mu}$-action on $M$ \cite{slicesubgroup}.  Define $$Y=G_{\mu}\times_K \left(\mathfrak{m}^*\times\left( N_1\oplus T_1\right)\right)$$  By the Symplectic Tube Theorem \ref{symplectic slice thm} a sufficiently small neighbourhood $Y_{\varepsilon}\subset Y$ of the zero section is equipped with a symplectic form $\omega_{Y_{\varepsilon}}$ and there is a $G_{\mu}$-invariant neighbourhood $U_0\subset M$ of $m$ and a $G_{\mu}$-equivariant symplectomorphism
\begin{equation*}
	\varphi: \left(Y_{\varepsilon},\omega_{Y_{\varepsilon}}\right) \longrightarrow \left( U_0,\restr{\omega}{U_0}\right)
\end{equation*}
with $\varphi\left([e,0]\right)=m$. We define 

$$\mathcal{N}=\left\lbrace (\rho,\; \nu+z_M(m))\in \mathfrak{m}^*\times (N_1\oplus T_1)\mid -\ad^*_z\mu+\rho\in\mathfrak{h}^{\circ}\right\rbrace.$$

By Proposition \ref{result A}, $\mathcal{N}$ is isomorphic to $N_1\oplus\mathcal{M}$, a $K$-vector space complementary to $\mathfrak{g}_{\mu}\cdot m$ in $\ker\left(D\J_H(m)\right)$.  Let $\Ne\subset \mathcal{N}$ be a $K$-invariant neighbourhood of $0$ such that $G_{\mu}\times_K\Ne\subset Y_{\varepsilon}$. We thus define $U\subset U_0$ by
\begin{equation}\label{RE: U}
	U:=\varphi\left(G_{\mu}\times_K\Ne\right)
\end{equation}
which, by Proposition \ref{result A},  is a $G_{\mu}$-invariant neighbourhood of $m$ in $\Phi_H^{-1}(\alpha)$.

Let $h_{\lambda}\in C^{\infty}(M)^H$ be an $H$-perturbation of $h$ with augmented Hamiltonian $h^{\xi}_{\lambda}$. We search for critical points of $h^{\xi}_{\lambda}$ in $U$. Using the local model \eqref{RE: U} this reads
$$h^{\xi}_{\lambda}:G_{\mu}\times_K\Ne\to \R.$$ As in the proof of Theorem \ref{crit} the critical points of $h^{\xi}_{\lambda}$ in $G_{\mu}\times_K\Ne$ lift to those of 
$$\widetilde{h}^{\xi}_{\lambda}:=\rho^*h^{\xi}_{\lambda}:G_{\mu}\times\Ne\to \R$$
where $\rho:G_{\mu}\times \Ne\to G_{\mu}\times_K\Ne$ is the orbit map. We may thus work with $\widetilde{h}^{\xi}_{\lambda}$ instead of $h^{\xi}_{\lambda}$. We define as usual the (left) action of the direct product $G_{\mu}\times K$ on $G_{\mu}\times \Ne$ by $$(h,k)\cdot (g,\nu)=(hgk^{-1},k\cdot \nu).$$ As $G_{\mu}\subset G_{\xi}$ by the \eqref{assumption} assumption, the lift $\widetilde{h}^{\xi}$ is $G_{\mu}\times K$-invariant whereas $\widetilde{h}^{\xi}_{\lambda}$ is only $H_{\mu}\times K$-invariant. 
By $\alpha$-nondegeneracy of $(m,\xi)$ and because $\varphi\left([e,0]\right)=m$:
\begin{equation*}\label{hknot2}
	\d \widetilde{h}^{\xi}(e,0)=0\quad\mbox{and}\quad D^2_{\mathcal{N}}\widetilde{h}^{\xi}(e,0)\quad\mbox{is non-singular}.
\end{equation*}
We use the Implicit Function Theorem and the compactness of $H_{\mu}\backslash G_{\mu}$ to get an $H_{\mu}$-invariant smooth function $\phi^{\eta}_{\lambda}:G_{\mu}\to \Ne$, depending on parameters $(\lambda, \eta)$ taken in a neighbourhood $V\subset \R\times \mathfrak{h}$ of $(0,\xi)$, satisfying
\begin{equation*}
		\d_{\mathcal{N}}\widetilde{h}^{\eta}_{\lambda}(g,\phi^{\eta}_{\lambda}(g))=0\quad\mbox{for every}\quad g\in G_{\mu}.
\end{equation*}

For every fixed parameters $(\lambda,\eta)\in V$, the $H_{\mu}\times K$-invariance of $\widetilde{h}^{\eta}_{\lambda}$ allows us to define a family of $H_{\mu}\times K$-invariant functions $f^{\eta}_{\lambda}$ by $$f^{\eta}_{\lambda}(g):=\widetilde{h}^{\eta}_{\lambda}(g,\phi^{\eta}_{\lambda}(g)).$$ Hence the implicit function $f^{\eta}_{\lambda}$ has a critical point at $g$ if and only if $\widetilde{h}^{\eta}_{\lambda}$ has a critical point at $(g,\phi^{\eta}_{\lambda}(g))\in G\times \Ne$. Then for $(\lambda,\eta)\in V$ the critical points of $h^{\eta}_{\lambda}$ in $U$ are in one-to-one correspondence with those of the function $f^{\eta}_{\lambda}$.
    \end{Proof}

\begin{cor}[Persistence of relative equilibria]\label{persistence of RE cor}
Under the assumptions of Theorem \ref{persistence of RE} the number of $H_{\mu}$-orbits of relative equilibria of $h$ with velocity close to $\xi$, that persist under a small $H$-perturbation in a neighbourhood of $G_{\mu}\cdot m$ in $\J_H^{-1}(\alpha)$, is bounded below by $\mbox{Cat}_{H_{\mu}}(G_{\mu}/G_m)$.
\end{cor}
\begin{Proof}
We apply Theorem \ref{gcrit} to $f^{\eta}_{\lambda}\in C^{\infty}(G_{\mu}/G_m)^{H_{\mu}}$ and we obtain that the number of $H_{\mu}$-orbits of critical points of $f^{\eta}_{\lambda}$ is bounded below by $\cat_{H_{\mu}}\left(G_{\mu}/G_m\right)$. In other words, as long as $(\lambda,\eta)\in V$, the number of $H_{\mu}$-orbits of relative equilibria of $h_{\lambda}$ with velocity $\eta$ in a neighbourhood of $G_{\mu}\cdot m$ in $\J_H^{-1}(\alpha)$ is at least $\cat_{H_{\mu}}\left(G_{\mu}/G_m\right)$.
\end{Proof}

\begin{example}[Torus action]\normalfont As a first application, we recover the result of \cite{Grabsi} for compact abelian groups and free actions. Let $(M,\omega,\T^n,\J_{\mbox{\tiny T}^n})$ be a hamiltonian $\T^n$-manifold where $\T^n$ is an $n$-dimensional torus acting freely on $M$ and let $\T^r$ be a subtorus of $\T^n$. Assume $h\in C^{\infty}(M)^{\mbox{\tiny T}^n}$ has an $\alpha$-nondegenerate relative equilibrium $(m,\xi)\in M\times \mathfrak{t}^r$ with momentum $\mu=\J_{\mbox{\tiny T}^n}(m)$ and where $\alpha=\restr{\mu}{\mathfrak{t}^r}$. As $\T^n$ and $\T^r$ are abelian, condition \eqref{assumption} always hold. Hence any $\T^r$-perturbation $h_{\lambda}$ with $\lambda$ small enough has at least $\cat_{\mbox{\tiny T}^r}(\T^n)$ $\T^r$-orbit of relative equilibria with velocity closed to $\xi$ in a neighbourhood of $\T^n\cdot m$ in $\J_{\mbox{\tiny T}^r}^{-1}(\alpha)$. Since $\T^n$ acts freely on $\T^r$ by left multiplication, $$\cat_{\mbox{\tiny T}^r}(\T^n)=\cat(\T^n/\T^r)=\cat(\T^{n-r}).$$ 
Hence $\cat_{\mbox{\tiny T}^r}(\T^n)=(n-r)+1.$
\end{example}

\begin{example}[Spherical pendulum on $S^3$]\label{s: RE examples}

As an application of Corollary \ref{persistence of RE cor}, we consider the case of a spherical pendulum on $S^3$, whose Hamiltonian is viewed as a perturbation of the free Hamiltonian on $S^3$. Endow $\R^4$ with the standard inner product $\langle\cdot,\cdot\rangle$ and let $e_1,e_2,e_3,e_4$ be the standard basis. The phase space is $(T^*S^3,\omega,G,\J_G)$ where $G=\SO(4)$ acts on 
\begin{equation*}
	T^*S^3=\left\lbrace (x,y)\in S^3\times \R^4\mid \langle x,y\rangle=0\right\rbrace
\end{equation*} 
by matrix multiplication $A\cdot (x,y)=(Ax,Ay)$. The associated momentum map $\J_G:T^*S^3\to \bigwedge^2(\R^4)$ is given by $\J_G(x,y)=y\wedge x.$ 

Let $H=\SO(3)\subset \SO(4)$ be the rotations about the $e_4$-axis with Lie algebra $\mathfrak{h}=\mathfrak{so}(3)_r$ as defined in Example \ref{ex:so(3) in so(4)}. The Hamiltonian of the spherical pendulum
\begin{equation}\label{hamiltonian pendulum}
	h_{\lambda}(x,y)=\frac{1}{2}\|y\|^2+\lambda \langle x,e_4\rangle
\end{equation}
 is an $H$-perturbation of the free Hamiltonian $h(x,y)=\frac{1}{2}\|y\|^2$. By definition, the relative equilibria of \eqref{hamiltonian pendulum} are pairs $\left((x,y),\xi\right)\in T^*S^3\times\mathfrak{h}$ such that
\begin{equation}\label{spherical pendulum RE equation}
		\d h_{\lambda}(x,y)=\d \phi_H^{\xi}(x,y)
\end{equation}
where $$\phi_{H}^{\xi}(x,y):=-\frac{1}{2}\mbox{Tr}((y\wedge x)\xi)=\langle \xi x,y\rangle.$$ Solving \eqref{spherical pendulum RE equation} is a straightforward calculation. The result is summarized as follows:
\begin{lemma}\label{lemma: RE pendulum}
	 With the notation of Example \ref{ex:so(3) in so(4)} we fix $\xi=(w,0)\in \mathfrak{h}$ for some $w\in\R^3$ and denote by $\pv(x)$ the projection of $x$ on $V=\mbox{span}(e_1,e_2,e_3)$. The relative equilibria $\left((x,y),\xi\right)\in T^*S^3\times\mathfrak{h}$ of \eqref{hamiltonian pendulum} satisfy the equations
	 \begin{enumerate}[label=(\roman*)]
	\item $\langle x,e_4\rangle=-\lambda\|w\|^{-2}$
	\item $\|\pv(x)\|^2=1-\langle x,e_4\rangle^2$ and $w\cdot \pv(x)=0$ (dot product in $\R^3$)
	\item $y=\xi x$ satisfies $\pv(y)=w\times \pv(x)$ and $\langle y,e_4\rangle=0$.
	\end{enumerate}
\end{lemma}

When $\lambda=0$ \eqref{hamiltonian pendulum} is the free Hamiltonian $h(x,y)=\frac{1}{2}\|y\|^2$ on $T^*S^3$. The integral curves of the corresponding hamiltonian vector field project to the great circles on $S^3$. We fix $\xi=(w,0)\in\mathfrak{h}$ with $w=(0,0,1)^T$. The relative equilibria $((x,y),\xi)$ of $h$ are such that $\langle x,e_4\rangle=0$ and $\pv(x)$ lies on the unit sphere in the hyperplane orthogonal to the line $[w]$, and $y$ is tangent to this sphere. In particular the pair $(m,\xi)$ with $m=(x,y)=(e_1,e_2)$ is a relative equilibrium of $h$. Its momentum is $\mu=\Phi_G(e_1, e_2)=\xi$ and its projection on $\mathfrak{h}^*$ is $\alpha=w^T=\begin{pmatrix}
0 & 0& 1\end{pmatrix}$. The stabilizer $G_{\mu}$ is a copy of $\SO(3)$ in $\SO(4)$ and the orbit $G_{\mu}\cdot m$ is the unit sphere $S^2\subset S^3$ lying on the hyperplane of equation $\langle x,e_4\rangle=0$.

We want to find the relative equilibria $((\xt,\yt),\eta)$ of the perturbed Hamiltonian \eqref{hamiltonian pendulum} which lie on $\J_H^{-1}(\alpha)$ where 
\begin{equation*}
	\J_H(\xt,\yt)=(\pv(\xt)\times\pv(\yt))^T
\end{equation*}
is the induced momentum map.  Writing $\eta=(u,0)\in \mathfrak{h}$ for some $u\in\R^3$, those relative equilibria satisfy the equation
\begin{equation}\label{pendulum: eq1}
	\|\pv(\xt)\|^2u-(\pv(\xt)\cdot u)\pv(\xt)=w
\end{equation}
with $w=(0,0,1)^T$, as fixed earlier. In addition they satisfy the equations of Lemma \ref{lemma: RE pendulum}, which require in particular that $\pv(\xt)\cdot u=0$. Replacing in \eqref{pendulum: eq1} we obtain $u=\|\pv(\xt)\|^{-2}w$. From Lemma \ref{lemma: RE pendulum} we get $\langle x,e_4\rangle=-\lambda\|\pv(\xt)\|^{4}$ and
\begin{equation}\label{eq:curve}
	\|\pv(\xt)\|^2+\lambda^2\|\pv(\xt)\|^8=1
\end{equation}
Setting $t=\|\pv(\xt)\|^2$ in \eqref{eq:curve}, we obtain the equation of an algebraic curve 
\begin{equation*}
	\lambda^2t^4+t-1=0\qquad t>0.
\end{equation*} 
For a fixed $\lambda$ there is exactly one solution representing the square of the radius $r(\lambda)$ of the sphere on which $\pv(\xt)$ lies. This sphere is an $H_{\mu}$-orbit of relative equilibria of \eqref{hamiltonian pendulum}. Since $r(0)=1$, it lies in a neighbourhood of the orbit $G_{\mu}\cdot m$ in $\J_H^{-1}(\alpha)$. Furthermore $\eta$ is such that $u=r(\lambda)^{-2}w$ and thus $\eta$ is close to $\xi$ in $\mathfrak{h}$. We also see from \eqref{eq:curve} that $\lambda$ must be chosen small enough such that $$\lambda<\|u\|^2<r(\lambda)^{-4}<c.$$ where $c$ is some constant coming from the fact that $r(\lambda)$ is bounded below.

We conclude that for $\lambda$ small enough, $h_{\lambda}$ has exactly one $H_{\mu}$-orbit of relative equilibria in a neighbourhood of $G_{\mu}\cdot m$ in $\J_H^{-1}(\alpha)$ with velocity close to $\xi$. For this example, we verify the assumptions of Theorem \ref{persistence of RE}. We have $G_{\mu}=H=SO(3)$ and the stabilizer $G_m$ is an $SO(2)$, as it is the subgroup of rotations in $SO(4)$ which preserve both axis $e_1$ and $e_2$. The quotient $G_{\mu}/G_m$ is thus a unit sphere $S^2$ and $H_{\mu}=G_{\mu}\cap H=SO(3)$. Furthermore, as $\mu=\xi$, the assumption \eqref{assumption} is satisfied, as well as the other assumptions of Theorem \ref{persistence of RE}. As expected, we have
\begin{equation*}
		\cat_{H_{\mu}}(G_{\mu}/G_m)=\cat_{SO(3)}(S^2)=1.
\end{equation*}
\end{example}

\section{Proof of Lemma \ref{hessian degeneracy}}\label{section hessian}
This section is devoted to the proof of Lemma \ref{hessian degeneracy} where we assume that $\mathfrak{g}$ is semi-simple. Further details are available in \cite{moi}. For each $\xi\in \mathfrak{g}$ a momentum map $\J_G:M\to \mathfrak{g}^*$ defines a smooth function $\phi_G^{\xi}:M\to\R$ depending linearly on $\xi$ $$\phi_G^{\xi}(m):=\langle \J_G(m),\xi\rangle.$$ Assume that $(m,\xi)\in M\times \mathfrak{g}$ is a relative equilibrium of some Hamiltonian $h\in C^{\infty}(M)^G$ with momentum $\mu=\J_G(m)$. By definition of a relative equilibrium, $\xi$ and $\mu$ commute i.e. $\ad^*_{\xi}\mu=0$. We would like to describe the space of degeneracy of the Hessian $D^2\phi_G^{\xi}(m)$ along the orbit $\mathfrak{g}\cdot m$. A straightforward calculation yields 
\begin{equation}
	D^2\phi_G^{\xi}(m)\left(y_M(m),x_M(m)\right) =\langle \J_G(m),[x,[y,\xi]]\rangle.
\end{equation}

Set $\mu=\J_G(m)$ and note that the Jacobi identity of the Lie bracket and the fact that $\ad^*_{\xi}\mu=0$ imply that $\langle \mu,[x,[y,\xi]]\rangle=\langle \mu,[y,[x,\xi]]\rangle$, reflecting the symmetric property of the Hessian. The degeneracy space of $D^2\phi_G^{\xi}(m)$ along $\mathfrak{g}\cdot m$ consists of the elements $y\in\mathfrak{g}$ such that
\begin{equation}\label{solve}
		\langle \mu,[y,[x,\xi]]\rangle=0\quad\mbox{for all}\quad x\in \mathfrak{g}.
\end{equation}
Since $\mu$ and $\xi$ commute, we can fix a maximal commutative Lie algebra $\mathfrak{t}\subset\mathfrak{g}$ such that $\xi\in\mathfrak{t}$ and $\mu\in\mathfrak{t}^*$. We complexify both of them
\begin{equation*}
	\mathfrak{g}_{\C}=\C\otimes_{\mathbb{R}}\mathfrak{g}\quad\mbox{and}\quad\mathfrak{t}_{\C}=\C\otimes_{\mathbb{R}}\mathfrak{t}
\end{equation*}
with extended Lie bracket $[\cdot,\cdot]_{\C}$. After this step the velocity and momentum read $\xi=1\otimes_{\mathbb{R}}\xi$ and $\mu=1\otimes_{\mathbb{R}}\mu$ and there respective stabilizer subalgebras are
\begin{equation*}
	\mathfrak{g}_{\xi}:=\left\lbrace x\in \mathfrak{g}_{\C}\mid [x,\xi]_{\C}=0\right\rbrace\quad\mbox{and}\quad\mathfrak{g}_{\mu}:=\left\lbrace x\in\mathfrak{g}_{\C}\mid \ad^*_x\mu=0\right\rbrace.
\end{equation*}
Consider the Cartan Lie subalgebra $\mathfrak{h}=\mathfrak{t}_{\C}$. Since $\xi\in\mathfrak{h}$ and $\mu\in\mathfrak{h}^*$, it is clear that $\mathfrak{h}$ is a subspace of both $\mathfrak{g}_{\xi}$ and $\mathfrak{g}_{\mu}$. We thus write
\begin{equation}\label{mu xi non regular}
	\mathfrak{g}_{\xi}=\mathfrak{h}\oplus\left(\bigoplus_{\beta\in S_f}\mathfrak{g}_{\beta}\right)\quad\mbox{and}\quad\mathfrak{g}_{\mu}=\mathfrak{h}\oplus\left(\bigoplus_{\alpha\in D_f}\mathfrak{g}_{\alpha}\right)
\end{equation}
for some finite subsets $S_f$ and $D_f$ of the root space $\mathcal{R}$ with the property: 
\begin{equation*}
	\alpha\in S_f\;(\mbox{resp. }D_f)\quad\Longrightarrow \quad-\alpha\in S_f \;(\mbox{resp. }D_f).
\end{equation*} 
\begin{definition}\label{regular}
	\normalfont $\xi$ (resp. $\mu$) is \defn{regular} if $S_f=\varnothing$ (resp. $D_f=\varnothing$). 
\end{definition}
Since $\mathfrak{g}_{\C}$ is semi-simple, the Killing form $\kappa$ induces an isomorphism $\kappa^{\sharp}:\mathfrak{h}^*\to \mathfrak{h}.$ Let $t_{\mu}\in \mathfrak{h}$ be the image of $\mu$ by this isomorphism and let $\mathcal{O}_{t_{\mu}}$ be the adjoint orbit of $t_{\mu}$. There is an identification
$$T_{t_{\mu}}\mathcal{O}_{t_{\mu}}=\sum_{\alpha\in\mathcal{R}\setminus D_f}\mathfrak{g}_{\alpha}.$$ The problem stated in \eqref{solve}, after complexification of the Lie algebra $\mathfrak{g}$, reduces to find all the $y\in \mathfrak{g}_{\C}$ satisfying 
\begin{equation}\label{non degeneracy with kappa}
		\kappa\left([y^*,t_{\mu}]_{\C},[x,\xi]_{\C}\right)=0\quad\mbox{for all}\quad x\in\mathfrak{g}_{\C}.
\end{equation}
Let $\left\lbrace H_1,\dots,H_k\right\rbrace\cup \left\lbrace X_{\alpha}\mid \alpha\in\mathcal{R}\right\rbrace$ be a Weyl-Chevalley basis of $\mathfrak{g}_{\C}$, where the $H_i$'s form a basis of $\mathfrak{h}$. Let $y\in\mathfrak{g}_{\C}$ be an arbitrary element and let $y^*=-\bar{y}$. With respect to the Weyl-Chevalley basis, this element is expressed as
\begin{equation}\label{y}
	y^*=\sum_{i=1}^ka_iH_i+\sum_{\alpha\in\mathcal{R}}\mu_{\alpha}X_{\alpha}\quad\mbox{for some unique}\quad a_i,\mu_{\alpha}\in\C.
\end{equation}
Hence
\begin{eqnarray*}
	[y^*,t_{\mu}]_{\C}&=&[\sum_{i=1}^ka_iH_i+\sum_{\alpha\in\mathcal{R}}\mu_{\alpha}X_{\alpha},t_{\mu}]_{\C}\\
	&=&\sum_{\alpha\in\mathcal{R}}\mu_{\alpha}[X_{\alpha},t_{\mu}]_{\C}\;\quad\mbox{as}\quad t_{\mu}\in\mathfrak{h}\\
	&=&-\sum_{\alpha\in\mathcal{R}}\mu_{\alpha}\alpha(t_{\mu})X_{\alpha}\\
	&=&-\sum_{\alpha\in\mathcal{R}\setminus D_f}\mu_{\alpha}\alpha(t_{\mu})X_{\alpha}
\end{eqnarray*}
where the last equality follows because $$[y^*,t_{\mu}]_{\C}\in T_{t_{\mu}}\mathcal{O}_{t_{\mu}}.$$ Similarly \eqref{mu xi non regular} allows us to write an element $[x,\xi]_{\C}\in T_{\xi}\mathcal{O}_{\xi}$ as $$[x,\xi]_{\C}=\sum_{\beta\in\mathcal{R}\setminus S_f}\lambda_{\beta}X_{\beta}\quad\mbox{with}\quad\lambda_{\beta}\in\C.$$ Solving \eqref{non degeneracy with kappa} is equivalent to solve
\begin{equation*}
		\sum_{\alpha\in\mathcal{R}\setminus D_f}\sum_{\beta\in\mathcal{R}\setminus S_f}\mu_{\alpha}\lambda_{\beta}\alpha(t_{\mu})\kappa\left(X_{\alpha},X_{\beta}\right)=0\quad\mbox{for any}\quad \lambda_{\beta}\in\C.
\end{equation*}
Using the fact that the $\mathfrak{g}_{\alpha}$'s appearing in the root decomposition are mutually orthogonal with respect to $\kappa$ (except for those corresponding to the same root with opposite sign), we get
\begin{eqnarray*}
		0&=&\sum_{\alpha\in\mathcal{R}\setminus D_f}\sum_{\beta\in\mathcal{R}\setminus S_f}\mu_{\alpha}\lambda_{\beta}\alpha(t_{\mu})\kappa\left(X_{\alpha},X_{\beta}\right)\\
		&=&\sum_{\alpha,\beta\in \mathcal{R}\setminus (D_f\cup S_f)}\mu_{\alpha}\lambda_{\beta}\alpha(t_{\mu})\kappa\left(X_{\alpha},X_{\beta}\right)\\
		&=&\sum_{\alpha\in \mathcal{R}\setminus (D_f\cup S_f)}\mu_{\alpha}\lambda_{\alpha}\alpha(t_{\mu})\kappa\left(X_{\alpha},X_{\alpha}\right)\\ &&+\sum_{\alpha\in \mathcal{R}\setminus (D_f\cup S_f)}\mu_{\alpha}\lambda_{-\alpha}\alpha(t_{\mu})\kappa\left(X_{\alpha},X_{-\alpha}\right)\\
		&=&\sum_{\alpha\in \mathcal{R}\setminus (D_f\cup S_f)}\mu_{\alpha}\alpha(t_{\mu})\left(\lambda_{\alpha}\kappa\left(X_{\alpha},X_{\alpha}\right)+\lambda_{-\alpha}\kappa\left(X_{\alpha},X_{-\alpha}\right)\right).\\
\end{eqnarray*}
This is true for any $\lambda_{\alpha}\in\C$ if and only if $\mu_{\alpha}=0$ for all $\alpha\in \mathcal{R}\setminus (D_f\cup S_f)$ as such roots satisfy $\alpha(t_{\mu})\neq 0$ and both $\kappa(X_{\alpha},X_{\alpha})$ and $\kappa(X_{\alpha},X_{-\alpha})$ do not vanish. We conclude that $y\in\mathfrak{g}_{\C}$ fulfils \eqref{non degeneracy with kappa} for all $x\in\mathfrak{g}_{\C}$ if and only if $y^*$ decomposes as
\begin{equation}
		y^*=\sum_{i=1}^ka_iH_i+\sum_{\alpha\in D_f\cup S_f}\mu_{\alpha}X_{\alpha}.
\end{equation}

\noindent Therefore, $$y^*\in \mathfrak{h}\oplus\left(\bigoplus_{\alpha\in D_f\cup S_f}\mathfrak{g}_{\alpha}\right)=\mathfrak{g}_{\xi}+\mathfrak{g}_{\mu}.$$
In particular this shows that the degeneracy set of the Hessian $D^2\Phi_G(m)$ along $\mathfrak{g}\cdot m$ belongs to $\mathfrak{g_{\xi}}+\mathfrak{g}_{\mu}$, by considering only the elements $y\in\mathfrak{g}_{\C}$ which are real. This proves the lemma because the other inclusion is clear.

\newcommand{\etalchar}[1]{$^{#1}$}

\vspace{0.5cm}
\begin{minipage}[t]{7cm}
MF:  {marine.fontaine@uantwerpen.be}\\
{\tt Departement Wiskunde-Informatica \\
Universiteit Antwerpen \\
2020 Antwerpen, BE.}

\end{minipage}\hfill
\begin{minipage}[t]{7cm}
JM:  {j.montaldi@manchester.ac.uk}\\
 {\tt School of Mathematics \\
University of Manchester \\
Manchester M13 9PL, UK.}
\end{minipage}
\end{document}